\documentclass[11pt]{article}
\usepackage{amssymb}
\font\tenBbb=msbm10
\def\Z{\hbox{\tenBbb Z}}
\def\T{{\hbox{\tenBbb T}}}
\font\tenBbb=msbm10

\def\sgn{\hbox{sgn}}
\font\tenBbb=msbm10
\def\N{\hbox{\tenBbb N}}
\font\tenBbb=msbm10
\newtheorem{theorem}{Theorem}[section]
\newtheorem{pro}{Proposition}[section]
\newtheorem{lem}{Lemma}[section]

\newtheorem{rem}{Remark}[section]
\newcommand{\R}{{ I\!\!R}}

\newcommand{\re}[1]{(\ref{#1})}

\def\into{\int \hspace*{-4mm} - \,}
\begin{document}
\begin{center}
\noindent {\Large \bf{Global well-posedness in $ L^2 $ for the
periodic Benjamin-Ono equation}}
\end{center}
\vskip0.2cm
\begin{center}
\noindent
{\bf Luc Molinet}\\
{\small L.A.G.A., Institut Galil\'ee, Universit\'e Paris-Nord,\\
93430 Villetaneuse, France.} \vskip0.3cm
E-mail : molinet@math.univ-paris13.fr.
\end{center}
\vskip0.5cm \noindent {\bf Abstract.} {\small We prove that the Benjamin-Ono equation is globally well-posed in $ H^s(\T) $
 for $ s\ge 0 $. Moreover we show that the associated flow-map is Lipschitz on every
  bounded set of $  H^s_0(\T) $, $s\ge 0$, and even real-analytic in this space
   for   small times.
  This result is sharp in the  sense that the flow-map (if it can be defined and coincides with the standard flow-map on $ {H}_0^\infty(\T) $)
   cannot be of class $ C^{1+\alpha} $, $\alpha>0 $, from $  H_0^s(\T) $ into $  H_0^s(\T) $
 as soon as  $ s< 0 $.}\vspace{2mm}\\

\section{Introduction, main results and notations}
\subsection{Introduction}
In this paper we continue  our  study (see \cite{M1}) of the Cauchy problem
for the
 Benjamin-Ono equation on the circle
$$
 \left\{\begin{array}{llll}
\partial_tu +{\cal{H}}\partial^2_x u - u \partial_x u=0 \,  \, ,\; (t,x)\in \R\times \T \;,\\
u(0,x)=u_0(x) \;,
\end{array}\right. \leqno{\mbox{(BO)}}
$$
where  $ \T=\R/2\pi \Z $, $ u $ is real-valued  and  ${\cal{H}}$ is the Hilbert transform defined for
  $ 2\pi  $-periodic functions with mean value zero by
$$
\widehat{{\cal{H}}(f)}(0):= 0 \quad \mbox{ and }\quad
\widehat{{\cal{H}}(f)}(\xi):= -i \,\sgn(\xi) \hat{f}(\xi), \quad
\xi\in \Z^* \quad \quad .
$$
\vskip0.3cm
\noindent
The Benjamin-Ono equation  arises as a model for long internal
gravity waves in deep stratified fluids, see \cite{B}.  This
equation possesses  a Lax pair  structure (cf. \cite{AF}, \cite{CW}) and thus
 has got
an infinite number of conservation laws. These conservation laws
permit to control  the $ H^{n/2} $-norms, $ n\in\N $, and thus to derive
 global well-posedness results in Sobolev spaces. The Cauchy
problem on the real line has been extensively studied these last
years (cf. \cite{Saut}, \cite{ABFS},  \cite{Io}, \cite{P},
\cite{MST}, \cite{KT1}, \cite{KK}). Recently, T. Tao \cite{T} has
 pushed the well-posedness theory to $ H^1(\R) $ by using an
appropriate gauge transform. This approach has been improved very
recently in \cite{PlaBu} and \cite{IoKe} where respectively $
H^s(\R) $, $s>0 $, and $ L^2(\R) $ are reached.  \\
 In the periodic setting, the local
well-posedness of (BO) is known in $ H^s(\T) $ for $s>3/2 $ (cf.
\cite{ABFS}, \cite{Io}), by standard compactness methods which do
not take advantage of the dispersive effects of the equation.
Thanks to the conservation laws  mentioned above and an
interpolation  argument, this leads to global well-posedness in $
H^{s}(\T) $ for $ s>3/2$ (cf. \cite{ABFS}). Very recently, F.
Ribaud and the author \cite{MR3} have improved the global
well-posedness result to $ H^1(\T) $ by using the gauge transform
introduced by T. Tao \cite{T} combining with Strichartz estimates
derived in \cite{B} for the Schr\"odinger group on the
one-dimensional torus. In \cite{M1} this approach combined with estimates in Bourgain type
 spaces leads to a global well-posedness result in the energy space $ H^{1/2}(\T) $.   Recall that
the Momentum and the Energy of  the Benjamin-Ono equation are respectively  given by

\begin{equation}
M(u):=\int_{\T}  u^2 \quad \mbox{ and } \quad
 E(u):=\frac{1}{2} \int_{\T} |D_x^{1/2} u |^2 +\frac{1}{6}
 \int_{\T}  u^3 \quad . \label{moment}
\end{equation}
The aim of this paper is to improve the local and global well-posedness
to $ L^2(\T) $.
\subsection{Notations}
For $x,y\in \R$,  $x\sim y$ means that there exists $C_1$,
$C_2>0$ such that $C_1 |x| \leq |y| \leq C_2|x|$ and $x\lesssim y$
means that there exists $C_2>0$ such that $ |x| \leq C_2|y|$. For
a Banach space $ X $, we denote by $ \| \cdot \|_X $ the norm in
$ X $.\\
We will use the same notations as in \cite{CKSTT1} and \cite{CKSTT2} to deal with Fourier transform of space periodic functions with a large period $\lambda $.
 $ (d\xi)_\lambda $ will be the renormalized counting measure on $ \lambda^{-1} \Z $ :
$$
\int a(\xi)\, (d\xi)_\lambda := \frac{1}{\lambda} \sum_{\xi\in
\lambda^{-1}\Z } a(\xi) \quad .
$$
As written in \cite{CKSTT2}, $ (d\xi)_\lambda $ is the counting measure on the integers when $ \lambda=1 $ and converges weakly
 to the Lebesgue measure when $ \lambda\to \infty $. In all the text, all the Lebesgue norms in $ \xi $ will be with respect to the
 measure $ (d\xi)_\lambda $. For  a $ (2\pi \lambda) $-periodic function $ \varphi$, we define its space Fourier transform on
  $ \lambda^{-1}\Z$ by
$$
\hat{\varphi}(\xi):=\int_{\R/(2\pi \lambda)\Z} e^{-i \xi x} \,
f(x) \, dx , \quad \forall \xi \in \lambda^{-1}\Z \quad .
$$
We denote by $ V(\cdot) $ the free group associated with the linearized Benjamin-Ono equation,
$$
\widehat{V(t)\varphi}(\xi):=e^{-i \xi|\xi|t} \,
\hat{\varphi}(\xi) , \quad \xi\in \lambda^{-1}\Z \quad .
$$
We define the Sobolev spaces $ H^s_{\lambda} $ for $ (2\pi\lambda)$-periodic functions by
$$
\|\varphi\|_{H^s_\lambda}:=\|\langle \xi \rangle^{s}
\widehat{\varphi}(\xi) \|_{L^2_\xi} =\|J^s_x \varphi
\|_{L^2_\lambda} \quad ,
$$
where $ \langle \cdot \rangle := (1+|\cdot|^2)^{1/2} $ and $
\widehat{J^s_x \varphi}(\xi):=\langle \xi \rangle^{s}
\widehat{\varphi}(\xi)
$. \\
For $ s\ge 0 $, the closed subspace of zero mean value functions of $ H^s_\lambda $ will be denoted by $ {H}^s_{0,\lambda} $ (it is equipped with the $ H^s_\lambda $-norm).\\
The Lebesgue spaces $ L^q_\lambda $, $ 1\le q \le \infty $, will be defined as usually by
$$
\|\varphi\|_{L^q_\lambda}:=\Bigl( \int_{\R/(2\pi\lambda) \Z}
|\varphi(x) |^q \, dx\Bigr)^{1/q}
$$
with the obvious modification for $ q=\infty $. \\
In the same way, for a function $ u(t,x) $ on $ \R\times\R/(2\pi \lambda) \Z $, we define its space-time Fourier transform by
$$
\hat{u}(\tau,\xi):={\cal F}_{t,x}(u)(\tau,\xi):=\int_{\R}
\int_{\R/(2\pi\lambda) \Z} e^{-i (\tau t+ \xi x)} \, u(t,x) \, dx
dt , \quad \forall (\tau,\xi) \in \R\times \lambda^{-1}\Z \quad .
$$
 We define the  Bourgain spaces $ X^{b,s}_\lambda $, $ Z^{b,s}_\lambda $,
  $ A_\lambda $  and $ Y^s_\lambda $ of $ (2\pi\lambda) $-periodic (in $x$) functions respectively
  endowed with the norm
\begin{eqnarray}
\| u \|_{X^{b,s}_{\lambda} } : =
 \| \langle \tau+\xi |\xi|\rangle^{b}  \langle \xi \rangle^s
  \hat{u}\|_{L^2_{\tau,\xi}} =
  \| \langle \tau\rangle^{b}  \langle \xi \rangle^s
  {\cal F}_{t,x}(V(-t) u ) \|_{L^2_{\tau,\xi}}\; ,
\end{eqnarray}
\begin{eqnarray}
\| u \|_{Z^{b,s}_{\lambda} } :=
 \|  \langle \tau+\xi |\xi|\rangle^b \langle \xi \rangle^s
  \hat{u} \|_{L^2_\xi L^1_\tau}= |  \langle \tau\rangle^b \langle \xi \rangle^s
  {\cal F}_{t,x}(V(-t) u ) \|_{L^2_\xi L^1_\tau}\; ,
\end{eqnarray}
\begin{eqnarray}
\| u \|_{A^b_\lambda}  & := &
 \| \langle \tau+\xi |\xi|\rangle^b \hat{u} \|_{L^1_{\tau,\xi}}=
  \| \langle \tau\rangle^b  {\cal F}_{t,x}(V(-t) u) \|_{L^1_{\tau,\xi}} \quad
\end{eqnarray}
and
\begin{equation}
\|u\|_{Y^s_\lambda} := \|u\|_{X^{1/2,s}_\lambda} +
\|u\|_{Z^{0,s}_\lambda} \quad ,
\end{equation}
where we will denote $ A^{0}_{\lambda} $ simply by $ A_{\lambda} $.
 Recall that  $ Y^s_{\lambda} \hookrightarrow Z^{0,s}_{\lambda}
 \hookrightarrow C(\R;H^s_\lambda) $. \\
We will also need the homogeneous semi-norm of $
\dot{X}^{b,s}_{\lambda} $ defined by
$$
\| u \|_{{\dot X}^{b,s}_{\lambda} }:  =
 \| |\tau+\xi |\xi ||^{b}  |\xi |^s
  \hat{u}\|_{L^2_{\tau,\xi}} \; .
$$
 $ L^p_t L^q_\lambda $  will denote the Lebesgue spaces
$$
\|u\|_{L^p_t L^q_\lambda}:=\Bigl( \int_{\R}
\|u(t,\cdot)\|_{L^q_\lambda}^p \, dt \Bigr)^{1/p}\quad
$$
with the obvious modification for $ p=\infty $.

Let $ u=\sum_{j\ge 0} \Delta_j u $ be a classical smooth non
homogeneous Littlewood-Paley decomposition in space of $ u $,
 $ \mbox{Supp } {\cal F}_x (\Delta_0 u ) \subset \R\times [-2,2] $
 and
$$
  \mbox{Supp }{\cal F}_x (\Delta_j u ) \subset \R\times [-2^{j+1},-2^{j-1}
 ]\cup \R\times  [2^{j-1},2^{j+1}] ), \quad j\ge 1\quad .
 $$
 We defined
the  Besov type space ${\tilde L}^4_{t,\lambda} $ by
 \begin{equation}
 \|u\|_{{\tilde L}^4_{t,\lambda}}:=
 \Bigl(\sum_{k\ge 0}
 \| \Delta_k u\|_{L^4_{t,\lambda}}^2\Bigr)^{1/2}
 \label{L4besov}
 \end{equation}
 Note that by the Littlewood-Paley square function theorem and Minkowski inequality,
 $$
 \|u\\|_{L^4_{t,\lambda}}\sim \Bigl\|\Bigl(\sum_{k=0}^\infty
 (\Delta_k
 u)^2\Bigr)^{1/2} \Bigr\|_{L^4_{t,\lambda}} \lesssim \Bigl(\sum_{k=0}^\infty
  \|\Delta_k u \|^2_{L^4_{t,\lambda}} \Bigr)^{1/2} =
  \|u\|_{{\tilde L}^4_{t,\lambda}}
  $$
  and thus $ {\tilde L}^4_{t,\lambda} \hookrightarrow
  L^4_{t,\lambda} $. \\
  We will work in the  function  spaces $ N_\lambda $ and $ {M}^s_\lambda $
   respectively defined  by
  $$
  \|u\|_{N_\lambda}:= \|u\|_{Z^{0,0}_\lambda}+
  \|Q_3 u\|_{X^{7/8,-1}_\lambda}+
 \|\chi_{[-4,4]}(t)\, u\|_{{\tilde L}^4_{t,\lambda}}
  $$
  and
  $$
  \|w\|_{M^s_\lambda}:=\|w\|_{Y^s_\lambda}+\|Q_1 w\|_{X^{1,-1}_\lambda} \quad ,
  $$
  where  $ Q_a $, $a\ge 0$, denotes the projection on the spatial Fourier modes of absolute value greater than $ a$ . \\
Finally,  for any function space $ B_\lambda $ and any $ T>0 $,
   we denote by  $ B_{T,\lambda} $ the corresponding restriction in time   space
 endowed with the norm
$$
\| u \|_{B_{T,\lambda}}: =\inf_{v\in B_{\lambda}} \{ \|
v\|_{B_{\lambda}} , \, v(\cdot)\equiv u(\cdot) \hbox{ on } ]0,T[
\, \}\;.
$$
 It is worth noticing that the map $ u\mapsto \overline{ u} $ is an isometry in all our function spaces.
 \\
We will denote by $ P_+ $ and $ P_- $ the projection on
respectiveley the positive and the negative spatial Fourier
modes. Moreover, for $ a\ge  0 $, we will denote by $ P_a $, $
Q_a $, $ P_{>a} $  and $ P_{<a} $ the projection on respectively
the spatial Fourier modes of absolute value equal or less than $
a$,
 the spatial Fourier modes of absolute value greater than $ a$, the  spatial
  Fourier modes  larger than $ a$ and the spatial Fourier modes smaller than $ a$.
\subsection{Main result}
Our well-posedness theorem reads :
\begin{theorem}\label{main}
For all $u_0\in H^{s}(\T)$ with $0\le s\le  1/2 $ and all $ T> 0 $, there
exists a  solution $ u $ of the Benjamin-Ono equation (BO) satisfying
\begin{equation}
u\in  C([0,T]; H^{s}(\T)) \cap N_{T,1}
 \quad \mbox{ and } \quad  \quad
P_+(e^{-i\partial_x^{-1} {\tilde u}/2} {\tilde u}) \in
X^{1/2,s}_{T,1}
 \label{class}
\end{equation}
 where
$$
 {\tilde u}:=u(t,x-t \into u_0)-\into u_0 \quad \mbox{ and  }\quad
 \widehat{\partial_x^{-1}}:= \frac{1}{i\xi}, \, \xi\in \Z^* \quad .
 $$
 This solution is unique in the class \re{class}. \\
Moreover $ u\in C_b(\R,L^{2}(\T)) $ and  the map $ u_0\mapsto u
$  is continuous from $H^{s}(\T)$ into $C([0,T],H^{s}(\T))$ and
Lipschitz on every bounded set from  $ {H}_0^{s}(\T) $ into
$C([0,T],{H}_0^{s}(\T))$.
\end{theorem}
Note that the result for $ s\ge 1/2 $ is established in \cite{M1}. Before stating our ill-posedness result let us make some comments  on Theorem \ref{main}.
\begin{rem}
We are not able to prove that  for any solution $ u $ of (BO) belonging  to
 $  C([0,T];  H^{s}(\T)) \cap N_{T,\lambda} $, the function
 $ P_+(e^{-i\partial_x^{-1} {\tilde u}/2} {\tilde u}) $ belongs to
  $  X^{1/2,s}_{T,\lambda} $. This is why we
  have to add this condition in our uniqueness class. Note however that
   any solution  that are limit  in $  C([0,T]; { H}^{s}(\T)) $
    of smooth
    solutions belongs to this class.
  Therefore, our solution  satisfies also the following (weaker) uniqueness notion
  used in \cite{IoKe} :
   it is the unique
   solution that is a  limit in $ C([0,T];H^s) $ of smooth solutions to (BO).
\end{rem}
\begin{rem}
Actually, we  prove that the flow-map is Lipschitz on every
bounded subset of  any hyperplan of $ H^{s}(\T) $ of functions with a fixed
mean value.
\end{rem}
\begin{rem}
The fact that $ u$ is real-valued is crucial to derive the
equation \re{eq2w} on $ w$. So, it does not seem that our approach
can be adapted to prove the local existence of complex-valued
solutions. On the other hand, it seems that a slight modification
of the  proof in \cite{M1} can lead to the local-wellposedness in
$ H^{1/2}(\T) $ for the  complex-valued version of (BO).
\end{rem}
Let us now state our ill-posedness issue.
\begin{theorem}\label{illposedness}
For $ s\ge 0 $ and $ t\in [0,1] $ the flow-map constructed by Theorem \ref{main} is
real-analytic from
 $ { H}_0^s(\T) $ into $ { H}_0^s(\T) $. On the other hand, for any $ t\in]0,1[ $ and any $ \alpha>0$,    the
 flow-map (if it can be defined and coincides with the standard flow-map on  $ { H}_0^\infty(\T)) $ cannot be of class $
 C^{1+\alpha} $ from $ { H}_0^s(\T) $ into ${
 H}_0^s(\T) $ as soon as $ s<0$.
 \end{theorem}
The main tools to prove Theorem \ref{main} are the gauge
transformation of T. Tao and the Fourier restriction spaces
introduced by Bourgain. Recall that in order to solve (BO), T.
Tao \cite{T} performed a kind of complex Cole-Hopf
transformation\footnote{Note that projecting (BO) on the non
negative frequencies, one gets the following equation : $
\partial_t (P_+ u) -i\partial_x^2 P_+ u =- P_+(u u_x) $} $
W=P_+(e^{-iF/2}) $, where $ F $ is a primitive of $ u $. In the
periodic setting, requiring that $ u $ has  mean value zero, we
can take  $ F=
  \partial_x^{-1} u $ the unique zero mean value primitive of $
  u$. By the mean value theorem, it is then easy to check
  that the above gauge transformation  is Lipschitz from $ L^2_\lambda $ to $
  L^\infty_\lambda $. This property, which is not true on the real line, is crucial to derive the
  smoothness of the flow-map. The equation satisfied by $ w =\partial_x W $
   takes the form
   $$
  w_t -i  w_{xx} = \partial_x P_+ (WP_- u_x )+ ...
   $$
   which  looks quite good since such nonlinear term enjoys   a strong smoothing effect on $ u$
     in Bourgain spaces.
   On the other hand, when one wants to inverse the gauge transformation, one gets something like
  $$
  u=e^{iF} w +...
  $$
  which is not so good since multiplication by  gauge function as
  $ e^{iF} $ behaves not so well in Bourgain spaces\footnote{Let us note that Bourgain spaces do not
   enjoy an algebra property}.
 Actually, the  ``bad'' regularity of $ u $ in the scale of Bourgain spaces
  is the main obstruction in going below $ H^{1/2}(\T) $ in \cite{M1}. In this work we substitute the above
    expression of $ u $ in
     the equation satisfied by $ w $. $ u $ still appears but only under the
      form $ e^{\mp iF/2} $ which possesses more regularities. On the other hand we
      have now to treat the multiplication by such functions in Bourgain spaces
       when estimating
       $ w$. Note that in the case $s=0 $ there is an additional difficulty
        mainly since we would like to control $ {\cal F}^{-1}_{t,x}(|\hat{u}|)
     $ in $ L^4_{t,x} $ whereas we only have a control on $ u $ in this space.
      This difficulty is overcome by noticing that actually $ u $ belongs
      to  a smaller space than $   L^4_{t,x} $ which is
    $       {\tilde L}^4_{t,x} $ (see \re{L4besov}).

   Concerning Theorem \ref{illposedness}, the fact that the
    flow-map (if it can be defined) cannot be of class $ C^3 $ in ${ H}_0^s(\T) $,
     $ s<0 $, can be
     obtained in the classical way for dispersive equations posed on $ \T $ (cf. \cite{Bo2}). To prove
      that it cannot be of class
     $ C^{1+\alpha} $, we somehow combine the bad behavior of the third iterate
       with the real-analyticity
     result in $ L^2(\T) $.
    \vspace*{5mm} \\
  {\it This paper is organized as follows:} In the next section
   we recall some linear estimates in Bourgain type spaces.
    In Section 3 we introduce the gauge transform and state the key nonlinear
     estimates. In Section 4, we prove the  estimates on the gauge function $ w$
      whereas the estimates on $ u$ are proven in Section \ref{Sec5}. In Section 4
    we derive uniform bounds for small
     initial data solutions and show a Lipschitz bound on the solution-map
      $ u_0\mapsto u $. The proof of Theorem \ref{main}  and
      Theorem \ref{illposedness} are completed respectively in
       Section 6 and  Section 7. Note that the proof of some technical
        lemmas needed in Sections \ref{Sec4}-\ref{Sec5} can be found in the appendix.
\section{Linear Estimates} \label{section2}
One of the main ingredient is the following linear estimate due to
Bourgain \cite{Bo1}.
\begin{equation}
\|v\|_{L^4(]-\pi,\pi[)L^4_1} \lesssim
\|v\|_{X^{3/8,0}_{]-\pi,\pi[,1}} \quad . \label{l5b}
\end{equation}
This estimate is proved in \cite{Bo1} (see also \cite{M1} for a
shorter proof) for Bourgain spaces of functions on $ \T^2 $ associated with the
Schr\"odinger group. The result for Bourgain space of functions on $ \R\times \T $
 can be  proven in exactly the same way (this can be  easily seen in the short proof presented in \cite{M1}).
 The corresponding estimate  for the Benjamin-Ono group follows by
writting $ v $ as the sum of its positive and negative spatial
modes parts. The estimate for any period $ \lambda\ge 1 $ follows
directly from dilation arguments. Indeed for any $ v\in
X^{3/8,0}_{1} $, setting $ v_\lambda :=\lambda^{-1} v(\lambda^{-2}
t,\lambda^{-1} x) $ , it is easy to see that $ v_\lambda\in
X^{3/8,0}_\lambda $ satisfies
$$
\|v_\lambda\|_{L^4_{t,\lambda}} =\lambda^{-1/4} \|v\|_{L^4_{t,1}}
\, , \quad \|v_\lambda\|_{{\dot
X}^{3/8,0}_{\lambda}}=\lambda^{-1/4} \|v\|_{{\dot
X}^{3/8,0}_{\lambda}}\mbox{ and }
\|v_\lambda\|_{L^2_{t,\lambda}}=\lambda^{1/2} \|v\|_{L^2_{t,1}}
\; .
$$
 From \re{l5b} we infer that for
 any function belonging to $  X^{3/8,0}_{\lambda} $ with $ \lambda\ge 1$,  it holds
\begin{equation}
\|v\|_{L^4_{t,\lambda}} \lesssim
\|v\|_{X^{3/8,0}_{\lambda}} \quad . \label{l5}
\end{equation}
 Let us now state some  estimates  for the free group
and the Duhamel operator. Let  $ \psi\in C_0^\infty([-2,2]) $ be
a time function  such that $ 0\le \psi \le 1 $ and $ \psi\equiv 1
$ on $ [-1,1] $. The following linear estimates are well-known (cf. \cite{Bo1}, \cite{G}).
\begin{lem}
\label{line2} For all  $ \varphi \in {H^s_\lambda}$ and all $ R>0
$, it holds :
\begin{equation}
\|\psi(t)V(t) \varphi \|_{Y^{s}_{\lambda}} \lesssim \|\varphi
\|_{H^s_\lambda} \quad , \label{l1}
\end{equation}
\begin{equation}
\|\psi(t/R) V(t) \varphi \|_{Z^{0,s}_{\lambda}} \lesssim \|\varphi
\|_{H^s_\lambda} \quad , \label{l2b}
\end{equation}
\begin{equation}
\|\psi(t/R) V(t) \varphi \|_{A_{\lambda}} \lesssim \|\hat{\varphi}
\|_{L^1_\xi} \quad , \label{l2}
\end{equation}
where it is worth noticing that the implicit constants in
\re{l2b} and \re{l2} do not depend on $ R $.
\end{lem}
{Proof. }\re{l1}  and \re{l2b} are classical. \re{l2} can be
obtained in the same way. Since $ V(t) $ commutes with any time
function and
$$ {\cal
F}_{x,t}(V(t)w(t,\cdot))=\hat{w}(\tau-\xi |\xi|,\xi)\quad ,
$$
  we infer that
\begin{eqnarray*}
 \|\psi(t/R) V(t)
\varphi \|_{A_{\lambda}}  & = & \|V(t)\psi(t/R) \varphi
\|_{A_{\lambda}}=\|{\cal F}_{t,x}(\psi(\cdot/R) \varphi
)\|_{L^1_{\tau,\xi}}
 \\
 & = & \| \hat{\psi}(\cdot)\|_{L^1_\tau}
 \|\hat{\varphi}\|_{L^1_\xi} \lesssim
 \|\hat{\varphi}\|_{L^1_\xi}\quad .
\end{eqnarray*}
Note that we  will use \re{l2b}-\re{l2} with $ R=\lambda^2 $ to estimate the low modes of $ u $  in  \re{eqbasfreq}.
\begin{lem}
\label{line3} For all $ G \in X^{-1/2,s}_{\lambda}\cap  Z^{-1,s}_{\lambda} $, it holds
\begin{equation}
\|\psi(t)\int_0^t V(t-t') G(t') \, dt' \|_{Y^{s}_{\lambda}}
\lesssim \| G \|_{X^{-1/2,s}_{\lambda}}+\| G
\|_{Z^{-1,s}_{\lambda}} \quad . \label{l3}
\end{equation}
\begin{equation}
\| \psi(t)\int_0^t V(t-t') G(t') \, dt' \|_{A_{\lambda}} \lesssim
\| G
 \|_{A^{-1}_{\lambda}} \quad . \label{l4b}
\end{equation}
\end{lem}
Let us recall that  \re{l3}-\re{l4b} are  direct consequences of
the following one dimensional (in time) inequalities (cf.
\cite{G} and \cite{GTV2}): for any function $ f\in {\cal S}(\R) $,
it holds
$$
\|\psi(t) \int_0^t f(t') \, dt' \|_{H^{1/2}_t} \lesssim
\|f\|_{H^{-1/2}_t} + \Bigl\|\frac{{\cal F}_t(f)}{\langle \tau
\rangle}\Bigr\|_{L^1_\tau} $$ and
$$
\Bigl\|{\cal F}_t \Bigl(\psi(t) \int_0^t f(t') \, dt'
\Bigr)\Bigr\|_{L^1_\tau} \lesssim  \Bigl\|\frac{{\cal
F}_t(f)}{\langle \tau \rangle}\Bigr\|_{L^1_\tau} \quad .
$$
\section{Gauge transform and nonlinear estimates} \label{Sec3}
\subsection{Gauge transform}
 Let $ \lambda\ge 1 $ and
 $ u $ be a smooth $(2\pi \lambda) $-periodic  solution of (BO) with initial data $ u_0 $.
 In the sequel,  we  assume that $u(t)  $ has  mean value zero for all time.
  Otherwise we do  the change of unknown :
 \begin{equation}
 v(t,x):=u(t,x-t \into u_0) -\into u_0 \label{chgtvar}\quad  ,
 \end{equation}
 where $ \into u_0:=P_0(u_0) =\frac{1}{2\pi \lambda} \int_{\R/(2\pi \lambda)\Z } u_0 $ is the mean value of $ u_0 $. It is easy to see that
 $ v $ satisfies (BO) with $ u_0-\into u_0 $ as initial data and since $ \into v $
  is preserved by the flow of (BO), $ v(t) $ has mean value zero for all time.
 We  define  $ F=\partial_x^{-1} u $ which is the periodic, zero mean value, primitive of $ u $,
 $$
 \hat {F}(0) =0 \quad \mbox{ and } \widehat{F}(\xi)=\frac{1}{ i\xi} \hat{u}(\xi) , \quad \xi\in \lambda^{-1}\Z^* \quad .
 $$
Following T. Tao \cite{T}, we introduce the gauge transform
\begin{equation}
W:=P_+(e^{-iF/2}) \quad . \label{defW}
\end{equation}
Since $ F $ satisfies $$ F_t +{\cal H} F_{xx}=\frac{F_x^2}{2}-\frac{1}{2} \into F_x^2=
  \frac{F_x^2}{2}-\frac{1}{2} P_0(F_x^2) \quad , $$
we  can  check that $ w:=W_x=-\frac{i}{2} P_+(e^{-iF/2} F_x)
=-\frac{i}{2} P_+(e^{-iF/2} u)$ satisfies
\begin{eqnarray}
w_t-iw_{xx} & = & -\partial_x P_+\Bigl[  e^{-iF/2}\Bigl(P_-(F_{xx})-\frac{i}{4} P_0(F_x^2)\Bigr)\Bigr] \nonumber \\
 & =  & -\partial_x P_+ \Bigl(W P_-( u_{x})  \Bigr)+ \frac{i}{4} P_0(F_x^2) w \label{eqw} \quad .
\end{eqnarray}
On the other hand, one can write $ u $ as
\begin{equation}
 u =  e^{iF/2} e^{-iF/2} F_x  = 2 i \, e^{iF/2}\partial_x (e^{-iF/2})= 2
i e^{iF/2} w  + 2 i  e^{iF/2} \partial_x P_- (e^{-iF/2} ) \quad . \label{A3}
\end{equation}
Recalling that $u $ is real-valued, we get
$$
u=\overline{u}=-2i e^{-iF/2} \overline{w} -2i e^{-iF/2} \partial_x \overline{P_-(e^{-iF/2})}
$$
and thus
\begin{equation}
P_-(u)=-2i P_-\Bigl(e^{-iF/2} \overline{w}\Bigr) -2i P_-\Bigl(e^{-iF/2} \partial_x P_+
(e^{iF/2})\Bigr)\quad
\label{eq2u}
\end{equation}
since $\overline{P_-(v)}=P_+(\overline{v})$ for any complex-valued
function $v $. Substituing \re{eq2u} in \re{eqw}, we obtain the
following equation satisfied by  $ w $ :
\begin{eqnarray}
w_t-iw_{xx} & =  & 2i\partial_x P_+ \Bigl(W \partial_x P_-(e^{-iF/2} \overline{w} )  \Bigr) \nonumber \\
& & \hspace*{-8mm}+ 2i\partial_x P_+ \Bigl[W \partial_x P_-\Bigl(
e^{-iF/2} \partial_x P_+ (e^{iF/2})\Bigr)\Bigr] +  \frac{i}{4}
P_0(F_x^2) W_x\label{eq2w} \quad .
\end{eqnarray}
Note also that it follows from \re{A3} that
\begin{eqnarray}
P_{>1} u & = & 2 i P_{>1} \Bigl( e^{iF/2}w\Bigr) + 2 i P_{>1} \Bigl( e^{iF/2} \partial_x P_-(e^{-iF/2}) \Bigr) \nonumber \\
& =  &  2 i P_{>1} \Bigl( e^{iF/2 }w\Bigr) + 2 i P_{>1} \Bigl(
P_{> 1}(e^{iF/2}) \partial_x  {P_-}(e^{-iF/2} ) \Bigr) \quad . \label{A4}
\end{eqnarray}
To end  this section is we state   the crucial nonlinear estimates on $ u $ and $ w $
 that will be proven in the next two sections.
  It is worth noticing that in all the estimates, we will replace the exponential function  (if it appears) by its entire serie and prove
 the absolute convergence of the resulting serie.
 Even if this approach can appear
  unecessary to prove the well-posedness result, it will be very useful
   in order to derive the analyticity of the flow-map. On the other hand it will
   in some estimates
    cause the appearance of a factor $ e^{\|
\widehat{\partial_x^{-1} u_0}\|_{L^1_\xi}} $ that could be avoid otherwise.  \\
\begin{pro} \label{NON}
Let $  u\in L^\infty_1 {H}^{s}_{0,\lambda}\cap N_{1,\lambda}$ be a
solution of (BO) and $  w\in X^{1/2,s}_{1,\lambda} $ satisfying
\re{eqw}-\re{A3}. Then for  $ 0\le s\le 1/2 $, it holds
\begin{eqnarray}
\|w\|_{M^{s}_{1,\lambda}} & \lesssim & (1+\|u_0\|_{L^2_\lambda})e^{\|
\widehat{\partial_x^{-1} u_0}\|_{L^1_\xi}}\;
\|u_0\|_{H^{s}_\lambda}  \nonumber \\
 & &  +
\|w\|_{X^{1/2,s}_{1,\lambda}} \Bigl(\|u\|_{N_{1,\lambda}}
+\|w\|_{X^{1/2,0}_{1,\lambda}} \Bigl) \, e^{\tilde K}
  \quad , \label{nonlinear1}
 \end{eqnarray}
\begin{equation}
\| u\|_{N_{1,\lambda}}  \lesssim  \|u_0\|_{L^2_\lambda}+ \Bigl(
\|w\|_{M^{0}_{1,\lambda}}+ \|u\|_{{N}_{1,\lambda}}^2\Bigr)e^{{\tilde K}}
 \quad
    \label{nonlinear2}
\end{equation}
and
\begin{equation}
\| u\|_{L^\infty_1 H^s_\lambda} \lesssim \|u_0\|_{L^2_\lambda}
+\Bigl(\|w\|_{M^{s}_{1,\lambda}} +  \|u\|_{N_{1,\lambda}}^2 \Bigr)
e^{{\tilde K}}
 \quad
\label{nonlinear4}
 \end{equation}
where
\begin{equation}
{\tilde K}=   C \, \Bigl(\|\widehat{\partial_x^{-1}
u_0}\|_{L^1_\xi}+ \|u\|_{N_{1,\lambda}}+\|u\|_{N_{1,\lambda}}^2
\Bigr)\quad \label{nonlinear3} \quad .
\end{equation}
for some universal constant $ C > 1 $. \\
\end{pro}
From Proposition \ref{NON} we will deduce uniform bounds for
smooth solutions of (BO) with small data (see Proposition
\ref{propuniformbound}). This will be the key point to derive the
local well-posedness result.
\section{Proof of the estimate on $ w$}  \label{Sec4}
In  this section, we will need the two following technical lemmas. The  first one, which is proven in the appendix, enables to  treat the multiplication with
 the gauge function $ e^{-iF/2} $
 in the Sobolev spaces  whereas
 the second  one (see the
  appendix of \cite{M1} for a proof), shows that, due
 to the frequency projections, we can share derivatives when taking the $ H^{s} $-norm of the second term of the right-hand side
 to \re{eq2w} or  \re{A4}.
\begin{lem}  \label{non1} Let $ 2\le q\le 4 $. Let $ h $ be  function of
$ H^1_\lambda $   and let $g \in L^q_\lambda $ such that $ J_x^\alpha g\in L^q_\lambda  $
 with $ 0< \alpha\le 1/2 $. Then it holds
\begin{equation}
\| J^\alpha (h g  ) \|_{L^{q}_\lambda} \lesssim \|J_x^\alpha
g\|_{L^q_\lambda} (\|h\|_{L^\infty_\lambda} +\|
h_x\|_{L^2_\lambda}) \quad . \label{non1a}
\end{equation}
\end{lem}
\begin{lem} \label{non2}
Let $ \alpha\ge 0 $  and $ 1<q<\infty $ then
\begin{equation}
\Bigl\|D^{\alpha}_x P_+\Bigl( f P_- \partial_x  g\Bigr) \Bigr\|_{L^q_\lambda} \lesssim
\| D^{\gamma_1}_x f \|_{L^{q_1}_\lambda } \,
  \| D^{\gamma_2}_x g \|_{L^{q_2}_\lambda} \quad , \label{estlemP+}
\end{equation}
with $ 1<  q_i<\infty $, $  1/q_1+1/q_2=1/q   $ and $
 \left\{ \begin{array}{l} \gamma_1 \ge \alpha , \; \gamma_2\ge 0 \\
 \gamma_1 + \gamma_2=\alpha+1 \end{array}\right. $ .
\end{lem}
\subsection{Choice of the extensions outside $ ]0,1[ $}
Let us introduce the following extensions outside the time interval  $ ]0,1[ $.
 Let $\tilde{\tilde w} $ be a zero-mean value extension of $ w $
 satisfying  $ \|\tilde{\tilde
 w}\|_{X^{1/2,0}_\lambda} \le 2 \|w\|_{X^{1/2,0}_{1,\lambda}} $ with
  $\tilde{\tilde w}=P_+( \tilde{\tilde w}) $,
  $ {\tilde W} $ be an extension of $ W $
 satisfying  $ \|{\tilde
 W}_x\|_{X^{1/2,s}_\lambda} \le 2 \|w\|_{X^{1/2,s}_{1,\lambda}} $ with
 $  {\tilde W}=P_+({\tilde W}) $ and let $ {\tilde w}:={\tilde W}_x $.
  We will also need a suitable extention $ {\tilde F } $ of $ F $.
  To construct $ \tilde{F} $ we proceed as follows : we take $ \tilde{\tilde u} $
  a zero-mean value extension of $ u $ in $ N_\lambda $ such that $ \|\tilde{\tilde u}\|_{N_\lambda}
  \le 2 \|u\|_{N_{1, \lambda}} $ and define $ \tilde{ u} $ by  setting
   $ Q_3 \tilde{u} = \psi Q_3 \tilde{\tilde u} $ and
   \begin{equation}
   \label{eqbasfreq}
   P_3 \tilde{u}=\psi(t/\lambda^2) \, P_3 V(t) u_0+\frac{\psi(t)}{2}  P_3 \Bigl[ \int_0^t V(t-t')
    \partial_x (\psi \tilde{\tilde u}(t'))^2 \, dt'\Bigr] \quad .
   \end{equation}
 The factor $ \lambda $ above will be very useful in \re{juju} to compensate  a factor $ \lambda $
  coming from the $ L^2_\lambda $-norm of $ \partial_x^{-1} u_0 $.
It is clear that  $ \tilde{u}  \equiv u $ on $ [0,1] $ and  $ \into \tilde{u}=0 $ on $\R $
 and thus we can  set $\tilde{F} = \partial_x^{-1} \tilde{u} $.

By the Duhamel formulation of \re{eq2w}, for $
 0\le t \le 1 $, we have
\begin{eqnarray}
w(t)& = & \psi(t) \Bigl[ V(t)w(0) +2i\int_0^t V(t-t')
\partial_x P_+ \Bigl((\psi\tilde{W}) \partial_x P_-(e^{-i\tilde{F}/2}
 \psi \overline{\tilde{\tilde w}} )  \Bigr) \nonumber \\
& & \hspace*{-8mm}+
2i   \int_0^t V(t-t')\partial_x P_+ \Bigl[ (\psi \tilde{W}) \partial_x P_-\Bigl(  e^{-i\tilde{F}/2} \partial_x P_+ ( e^{i\tilde{F}/2})\Bigr)\Bigr]\nonumber \\
 & & \hspace*{1cm} +\frac{i}{4} \int_0^t V(t-t')\Bigl(
  P_0( \tilde{u}^2) \psi \tilde{W}_x\Bigr)(t')\, dt' \, \Bigr] \quad . \label{wetoile}
 \end{eqnarray}
 To obtain the desired estimates we will first apply Lemmas \ref{line2}-\ref{line3} to \re{wetoile} and  then apply Lemmas \ref{nonlin1}-\ref{nonlin3} below with
 $ W:=\psi {\tilde W} $, $ F:={\tilde F} $ and $ v:= \psi \overline{{\tilde {\tilde w}}} $. Note that since
  ${\tilde {\tilde w}}=P_+ {\tilde {\tilde w}} $,
   we have $\overline{{\tilde {\tilde w}}}=P_- \overline{{\tilde {\tilde w}}} $ and thus
    $ v=P_- v $. Moreover,
    $ W $ and $ v$ being supported in time in $ \{t\in \R, |t|\le 2\} $, $ W=\psi_2 W $ and $ v=\psi_2 v $ where $ \psi_2(\cdot)=\psi(\cdot/2)$ and $ \psi $  is the cut-off in time function
     defined in Section \ref{section2}.
\subsection{Some multilinear estimates}
 The main tool for proving \re{nonlinear1} are  three multilinear estimates.
These estimates enlight the good behavior in Bourgain spaces of the  terms
 of the right-hand side  of \re{eq2w}. In the following lemmas
  $W $, $w:=\partial_x W $  and $ v$   are assumed to be supported
  in time in $ [-2,2] $ and  we set $ \psi_2(\cdot)=\psi(\cdot/2)$ (see above).
\begin{lem} \label{nonlin1}
For any $ s\ge 0 $ and $ 0<\varepsilon <\!\! < 1 $,
\begin{equation}
\Bigl\|\partial_x P_+ \Bigl[W \partial_x P_-\Bigl(e^{-iF/2}
\partial_x P_+ (e^{iF/2})\Bigr)\Bigr] \Bigr\|_{X^{-1/2\,
+\varepsilon,s}_{\lambda} }\lesssim  \|w \|_{X^{1/2,s}_{\lambda}}
\|\psi_2  \, F_x \|_{L^4_{t,\lambda}}^2 e^{C\,
\|F\|_{L^\infty_{t,\lambda}}}\label{C1}\quad .
\end{equation}
\end{lem}
{\it Proof. }
 As written above, we will actually prove \re{C1} with as left-hand side member
 (Note that the factor $
 e^{\|F\|_{L^\infty_{t,\lambda}}} $ in the right-hand side of \re{C1} could be  avoid otherwise):
$$
\sum_{k\ge 1} \sum_{l\ge 1} \frac{1}{k!}\frac{1}{l!} \Bigl\|
\partial_x P_+ \Bigl[W \partial_x P_-\Bigl(F^k \partial_x P_+
 F^l\Bigr)\Bigr] \Bigr\|_{X^{-1/2+\varepsilon,s}_{\lambda}
}\quad .
$$
Note that, according to  the  support in time  of $ W $,  the expression  contained in norm
 remains unchanged by multiplication with
  the cut-off in time function $ \psi_2 $.
Setting
$$g=\partial_x P_-\Bigl(\psi_2  F^k\partial_x P_+ (\psi_2 F^l)\Bigr)\; ,
 $$ it follows from Lemma \ref{non2} that
$$
\|g\|_{L^2_{t,\lambda}} \lesssim
\|\psi_2\partial_x(F^k)\|_{L^4_{t,\lambda}}\|\psi_2\partial_x(F^l)\|_{L^4_{t,\lambda}}
\lesssim k \, l \|\psi_2 F_x\|_{L^4_{t,\lambda}}^2
\|F\|_{L^\infty_{t,\lambda}}^{k+l-2} \;.
$$
 It thus suffices   to prove that
\begin{equation}
\| \partial_x P_+(W P_- g )\|_{X^{-1/2 +\varepsilon,s }} \lesssim  \|
w\|_{X^{1/2,s}} \|g \|_{L^2_{t,\lambda}} \quad .
\label{mi1}
\end{equation}
 By duality it is equivalent to estimate
\begin{eqnarray*}
I & =& \Bigl| \int_{A}  \xi\,  \hat{h}(\tau,\xi)
  \xi_1^{-1} \hat{w}(\tau_1,\xi_1) \hat{g}(\tau_2,\xi_2) \Bigr|
\end{eqnarray*}
where $ (\tau_2,\xi_2)=(\tau-\tau_1,\xi-\xi_1) $,
and due to the frequency projections
$$
A=\{(\tau,\tau_1, \xi,\xi_1)\in \R^2\times(\lambda^{-1}\Z)^2,
\quad \xi\ge 1/\lambda, \, \xi_1\ge 1/\lambda , \, \xi-\xi_1\le
-1/\lambda \quad \} \quad.
$$
 Note that in the domain of integration above,
\begin{equation}
\xi_1 \ge |\xi-\xi_1| \quad \mbox{ and } \quad \xi_1\ge \xi \quad
. \label{C3}
\end{equation}
It thus folllows that
$$
I\lesssim  \int_{A}  \langle \xi\rangle^{-s} |\hat{h}(\tau,\xi)| \langle \xi_1\rangle^{s} |\hat{w}(\tau_1,\xi_1)| |\hat{g}(\tau_2,\xi_2)|
$$
 and on account of \re{l5},
\begin{eqnarray*}
I & \lesssim & \| {\cal F}^{-1} (\langle \xi\rangle^{-s}|\hat{h}| )\|_{L^4_{t,\lambda}}
\| {\cal F}^{-1} (\langle \xi\rangle^{s}|\hat{w}| )\|_{L^4_{t,\lambda}}
\|{\cal F}^{-1} (|\hat{g}| )\|_{L^2_{t,\lambda}} \\
 & \lesssim &  \,\|h\|_{X^{3/8,-s}_\lambda}  \|w\|_{X^{1/2,s}_\lambda} \|g \|_{L^2_{t,\lambda}}
\end{eqnarray*}
which proves \re{C1}.
\begin{lem} For any $ s\ge 0 $ it holds
\begin{eqnarray}
\hspace*{-15mm}\Bigl\|  \partial_x P_+ \Bigl( W\partial_x
P_-(e^{-iF/2} P_- v )\Bigr)\Bigr\|_{X^{-1/2,s}_\lambda}
 & \lesssim &   \|w\|_{X^{1/2,s}_\lambda} \|v\|_{X^{1/2,0}_\lambda} \, e^{C \,\|F\|_{A_\lambda}} \nonumber \\
& & \hspace*{-70mm}\Bigl( 1+\|P_3 F\|_{{\dot X}^{1,0}_\lambda} +\|
P_{>3} F_x \|_{X^{7/8,-1}_\lambda} + \|F\|_{A_\lambda}+ \|\psi_2
F_x \|_{\tilde{L^4}_{t,\lambda}} \Bigr)\; . \label{mi2}
\end{eqnarray}
\label{nonlin2}
 \end{lem}
{\it Proof. } Again we will in fact prove \re{mi2} with as
left-hand side member :
\begin{eqnarray*}
 \Bigl\|  \partial_x P_+ \Bigl( W\partial_x  P_- v )\Bigr)\Bigr\|_{X^{-1/2,s}_\lambda}
 +\sum_{k\ge 1}\ \frac{1}{k!}
 \Bigl\|  \partial_x P_+ \Bigl( W\partial_x P_-(F^k P_- v )\Bigr)\Bigr\|_{X^{-1/2,s}_\lambda} \quad.
\end{eqnarray*}
The first term of the above inequality is estimated in (\cite{M1},
Lemma 3.3) by
\begin{equation}
\Bigl\|  \partial_x P_+ \Bigl( W\partial_x  P_- v
)\Bigr)\Bigr\|_{X^{-1/2,s}_\lambda} \lesssim  \|w\|_{X^{1/2,s}}
\|v\|_{X^{1/2,0}_\lambda}\quad .\label{LL1}
\end{equation}
By duality it thus remains to estimate
\begin{eqnarray}
D_k & =& \Bigl| \int_{B}  \xi\,  \hat{h}(\tau,\xi)
  \xi_1^{-1} \hat{w}(\tau_1,\xi_1)(\xi-\xi_1)
  \widehat{P_- v}(\tau_2,\xi_2)\prod_{i=3}^{k+2}\hat{F}(\tau_i,\xi_i) \Bigr|\label{C2}
\end{eqnarray}
where $
(\tau_{k+2},\xi_{k+2})=(\tau,\xi)-\sum_{i=1}^{k+1}(\tau_i,\xi_i)
$, and due to the frequency projections
\begin{eqnarray*}
B=\{(\tau,\tau_1,..,\tau_{k+1}, \xi,\xi_1,..,\xi_{k+1})\in
\R^{k+2}\times(\lambda^{-1}\Z)^{k+2},
 & & \\
& & \hspace*{-60mm} \xi_1\ge \xi\ge 1/\lambda, \, \xi-\xi_1\le -1/\lambda
\quad \} \quad.
\end{eqnarray*}
First splitting $ D_k $ into the two following terms $$ I_k =
\Bigl| \int_{B}  \xi\,  \hat{h}(\tau,\xi)
  \xi_1^{-1} \hat{w}(\tau_1,\xi_1)(\xi-\xi_1)
  \widehat{P_{\{2^{10}k\}} P_- v}(\tau_2,\xi_2)\prod_{i=3}^{k+2}\hat{F}(\tau_i,\xi_i) \Bigr|
$$
and
$$J_k  = \Bigl| \int_{B}  \xi\,  \hat{h}(\tau,\xi)
  \xi_1^{-1} \hat{w}(\tau_1,\xi_1)(\xi-\xi_1)
  \widehat{Q_{\{2^{10}k\}} P_- v}(\tau_2,\xi_2)\prod_{i=3}^{k+2}\hat{F}(\tau_i,\xi_i) \Bigr|
  $$
 we observe that
 \arraycolsep1pt
\begin{eqnarray}
I_k &
  \lesssim &  \| {\cal F}^{-1} (\langle
\xi\rangle^{-s}|\hat{h}| )\|_{L^4_{t,\lambda}} \| {\cal F}^{-1}
(\langle \xi\rangle^{s}|\hat{w}| )\|_{L^4_{t,\lambda}}
\Bigl\|\partial_x \Bigl( (P_{\{2^{10}k\}} P_- v) F^k\Bigr)
\Bigr\|_{L^2_{t,\lambda}} \nonumber \\
&  \lesssim &   k\,
\|h\|_{X^{3/8,-s}_\lambda}\|w\|_{X^{1/2,s}_\lambda}(\|v\|_{L^2_{t,\lambda}}
\|F\|_{L^\infty_{t,\lambda}}+\|v\|_{L^4_{t,\lambda}}\|\psi_2
F_x\|_{L^4_{t,\lambda}})
 \|F\|_{L^\infty_{t,\lambda}}^{k-1} \label{mi8bbb}\quad ,
\end{eqnarray}
\arraycolsep3pt
 since obviously,
\begin{eqnarray*}
\Bigl\|\partial_x \Bigl((P_{\{2^{10}k\}} P_- v) F^k
\Bigr)\Bigr\|_{L^2_{t,\lambda}} & \lesssim & \|P_{\{2^{10}k\}} P_-
v_x\|_{L^2_{t,\lambda}}
 \|F\|_{L^\infty_{t,\lambda}}^k \\
  & & + k \|P_{\{2^{10}k\}} P_- v\|_{L^4_{t,\lambda}}
 \| \psi_2 F_x\|_{L^4_{t,\lambda}}
 \|F\|_{L^\infty_{t,\lambda}}^{k-1} \\
  &\lesssim &  k(\| v\|_{L^2_{t,\lambda}} \|F\|_{L^\infty_{t,\lambda}}
  +\|  v\|_{L^4_{t,\lambda}} \|\psi_2 F_x\|_{L^4_{t,\lambda}} )
 \|F\|_{L^\infty_{t,\lambda}}^{k-1} \quad .
\end{eqnarray*}
It thus remains  to estimate $ J_k $.
 Note  that since \re{C3} holds on $ B $,
  setting
  $$
B_1=\{(\tau,\tau_1,..,\tau_{k+1}, \xi,\xi_1,..,\xi_{k+1})\in B,\,
|\xi|\le 2^{10} k  \mbox{ or } |\xi-\xi_1|\le 2^{10} k \, \}
  $$ we
get thanks to \re{l5},
 \begin{eqnarray}
{J_{k_{/B_1}}}  & \lesssim & k\, \| {\cal F}^{-1} (\langle
\xi\rangle^{-s}|\hat{h}| )\|_{L^4_{t,\lambda}} \| {\cal F}^{-1}
(\langle \xi\rangle^{s}|\hat{w}| )\|_{L^4_{t,\lambda}} \|
(Q_{\{2^{10}k\}}P_-v) F^k \|_{L^2_{t,\lambda}}
 \nonumber\\
& \lesssim & k\,
\|h\|_{X^{3/8,-s}_\lambda}\|w\|_{X^{1/2,s}_\lambda}\|v\|_{L^2_{t,\lambda}}
 \|F\|_{L^\infty_{t,\lambda}}^k \label{mi8bb}\quad .
\end{eqnarray}
It thus suffices  to control
\begin{equation}
{ J_{k_{/B_2}}}  = \Bigl| \int_{B_2}  \xi\, \hat{h}(\tau,\xi)
  \xi_1^{-1} \hat{w}(\tau_1,\xi_1)(\xi-\xi_1)
  \widehat{(Q_{\{2^{10}k\}}P_- v)}(\tau_2,\xi_2)\prod_{i=3}^{k+2}\hat{F}(\tau_i,\xi_i) \Bigr|
  \label{Cc2}
\end{equation}
where
$$
 B_2=\{(\tau,\tau_1,..,\tau_{k+1}, \xi,\xi_1,..,\xi_{k+1})\in B,\,
\xi> 2^{10} k  , \, \xi-\xi_1<- 2^{10} k \, \}\; .
$$
 One of the main difficulties will be  that we do not
have a control on $ \|{\cal F}^{-1}_{t,x} (|\hat{F_x}|
)\|_{L^4_{t,\lambda}} $ but only on
 $ \|F_x\|_{L^4_{t,\lambda}} $. This can be  overcame when $ s>0 $ but causes
 a kind of logarithmic divergence when $s=0 $. To control $ J_{k_{/B_2}}  $
we will have to use the stronger norm $ \tilde{L}^4_{t,\lambda} $
of $ F_x $. To simplify the notation we denote $ Q_{\{2^{10}k\}}
P_- v $ by $ {\tilde v} $. Since   we cannot ``force"
 the integrant to be non negative in \re{Cc2}, we have to act
 carefully.
 We notice that using Littlewood-Paley decomposition (see \re{L4besov})
 we can rewrite $ Q_{2^{10}k}({\tilde v} F^k) $ as
 \begin{eqnarray*}
 Q_{\{2^{10}k\}}({\tilde v} F^k)& = &   Q_{\{2^{10}k\}}\Bigl( \sum_{i_2\ge 8+\alpha(k)}
 \hspace*{-2mm}\Delta_{i_2}({\tilde v})\hspace*{-2mm}
 \sum_{ i_3\ge i_2-6-\alpha(k)} \hspace*{-5mm}\Delta_{i_3} (F)\hspace*{-2mm}
  \sum_{0\le
 i_4,..,i_{k+2}\le i_3}\hspace*{-5mm} n(i_3,..,i_{k+2})\prod_{j=4}^{k+2}\Delta_{i_j} (F) \Bigr) \\
  & & +   Q_{\{2^{10}k\}}\Bigl( \sum_{i_2\ge 8+\alpha(k)}\hspace*{-2mm} \Delta_{i_2}({\tilde v})
 \hspace*{-4mm} \sum_{0\le i_3< i_2-6-\alpha(k)}
\hspace*{-5mm} \Delta_{i_3} (F) \hspace*{-2mm}\sum_{0\le
 i_4,..,i_{k+2}\le i_3}\hspace*{-5mm} n(i_3,..,i_{k+2})\prod_{j=4}^{k+2} \Delta_{i_j} (F)
 \Bigr)\\
 & =& G_1+G_2 \quad ,
 \end{eqnarray*}
 where $\alpha(k) $ denotes the entire part of  $\ln(k)/\ln(2) $ and
$n(i_3,..,i_{k+2})$ is an integer belonging to $ \{1,..,k\} $ (Note for instance that $ n(i_3,..,i_{k+2})=1 $
 for $ i_3=i_4=\cdot\cdot=i_{k+2} $ and $ n(i_3,..,i_{k+2})=k $ for
$ i_3\neq i_4\neq\cdot\cdot\neq i_{k+2} $).
  From \re{Cc2} we thus infer that
\begin{eqnarray}
{J_{k_{/B_2}}} & \lesssim & \sum_{i=1}^2 \int_{B_1}  \xi\,
|\hat{h}(\tau,\xi)|
  \xi_1^{-1} |\hat{w}(\tau_1,\xi_1| |\xi-\xi_1| |\widehat{G_i}(\tau-\tau_1,\xi-\xi_1)|
  \nonumber \\
 & \lesssim &  \Lambda_1+\Lambda_2\; . \label{dodo}
\end{eqnarray}
$ \bullet ${\it Estimate on $ \Lambda_1 $.}
  Thanks to  the definition of $ B $ and \re{l5}, we easily obtain
   \begin{eqnarray*}
\Lambda_1  & \lesssim &  \| {\cal F}^{-1} (\langle
\xi\rangle^{-s}|\hat{h}| )\|_{L^4_{t,\lambda}} \| {\cal F}^{-1}
(\langle \xi\rangle^{s}|\hat{w}| )\|_{L^4_{t,\lambda}} \|
\partial_x G_1 \|_{L^2_{t,\lambda}}
 \nonumber\\
& \lesssim &
\|h\|_{X^{3/8,-s}_\lambda}\|w\|_{X^{1/2,s}_\lambda}\|\partial_x G_1\|_{L^2_{t,\lambda}}
  \quad .
\end{eqnarray*}
  On the other hand, using the frequency support of the functions,  we infer that for $
q\ge 9+\alpha(k) $,
$$
\Delta_q(G_1)  =   Q_{\{2^{10}k\}}\Delta_q \Bigl( \sum_{i_3\ge
q-8-\alpha(k)\atop i_3\ge 2  } \hspace*{-5mm}\Delta_{i_3} (F)
 \sum_{ i_2\ge 8+\alpha(k)\atop i_2\le i_3+6+\alpha(k)}\hspace*{-5mm}  \Delta_{i_2}({\tilde v}) \hspace*{-5mm}\sum_{0\le
 i_4,..,i_{k+2}\le i_3} \hspace*{-5mm} n(i_3,..,i_{k+2})\prod_{j=4}^{k+2} \Delta_{i_j} (F)
 \Bigr)
$$
and thus
  $$
  \|\Delta_q G_1\|_{L^2_{t,\lambda}}\lesssim
  \sum_{i_3\ge q-8-\alpha(k) \atop i_3\ge 2} \|\psi_2 \, \Delta_{i_3} F\|_{L^4_{t,\lambda}}
  \|{\tilde v}
  \|_{L^4_{t,\lambda}}\Bigl\| \sum_{0\le
 i_4,..,i_{k+2}\le i_3} \hspace*{-5mm} n(i_3,..,i_{k+2})\prod_{j=4}^{k+2} \Delta_{i_j}
 (F)\Bigr\|_{L^\infty_{t,\lambda}}\; .
 $$
 But
\begin{eqnarray}
\Bigl\| \sum_{0\le
 i_4,..,i_{k+2}\le i_3} \hspace*{-5mm}n(i_3,..,i_{k+2})\prod_{j=4}^{k+2} \Delta_{i_j}
 (F)\Bigr\|_{L^\infty_{t,\lambda}} &\lesssim& k\, \Bigl\|
 \sum_{0\le
 i_4,..,i_{k+2}}  \hspace*{-5mm}|\widehat{\Delta_{i_4}
 (F)}|*..*  |\widehat{\Delta_{i_{k+2}}
 (F)}|\Bigr\|_{L^1_{\tau,\xi}} \nonumber \\
& \lesssim& k\, \Bigl\|\Bigl(
 \sum_{i_4\ge 0 }|\widehat{ \Delta_{i_4}
 (F)}| \Bigr) *..*\Bigl(
 \sum_{i_{k+2}\ge 0 } |\widehat{\Delta_{i_{k+2}}
 (F)}| \Bigr) \Bigr\|_{L^1_{\tau,\xi}}\nonumber \\
  & \lesssim & k\,
 \|F\|_{A_\lambda}^{k-1} \label{estA}\quad .
\end{eqnarray}
Therefore,
 \begin{eqnarray*}
\|\partial_x G_1\|_{L^2_{t,\lambda}}^2 &  \sim  &  \sum_{q\ge
9+\alpha(k)} 2^{2q}
 \|\Delta_q G_1\|^2_{L^2_{t,\lambda}} \\
 & \lesssim  & k^2 \|{\tilde v}\|_{L^4_{t,\lambda}}^2
 \|F\|^{2(k-1)}_{A_{\lambda}}
\sum_{q\ge 9+\alpha(k)} \Bigl( \sum_{j\ge q-8-\alpha(k)\atop j\ge
2} 2^{(q-j)} 2^j \|\psi_2 \Delta_j F\|_{L^4_{t,\lambda}} \Bigr)^2
\end{eqnarray*}
 But, by the definition of the norm $ {\tilde L}^4_{t,\lambda} $
 \footnote{Note that we could avoid the $ {\tilde L}^4_{t,\lambda} $-norm here by
  invoking the Littlewood-Paley square function theorem in the estimate on $ G_1 $} (see \re{L4besov}),
 for $ j\ge 2 $, $ 2^j\|\Delta_j F\|_{L^4_{t,\lambda}}\lesssim \gamma_j \|F_x\|_{{\tilde L}^4_{t,\lambda}} $
  with $ \|(\gamma_j)\|_{l^2(\N)}\lesssim 1 $. Hence, by Young inequality,
  $$
   \sum_{j\ge q-8-\alpha(k)\atop j\ge 2} 2^{(q-j)} 2^j \|\psi_2
   \Delta_j F\|_{L^4_{t,\lambda}}
   \lesssim k \gamma_q  \|\psi_2 F_x\|_{{\tilde L}^4_{t,\lambda}}
  $$
  and thus
 $$
 \|\partial_x G_1\|_{L^2_{t,\lambda}} \lesssim k^2
 \|{\tilde v}\|_{L^4_{t,\lambda}} \|\psi_2 F_x\|_{{\tilde
L}^4_{t,\lambda}}
 \|F\|^{k-1}_{A_{\lambda}}  \quad .
 $$
 Therefore,  the following estimate holds
 \begin{equation}
 \Lambda_1  \lesssim k^2
\|h\|_{X^{3/8,-s}_\lambda}\|w\|_{X^{1/2,s}_\lambda}
   \|v\|_{L^4_{t,\lambda}} \|\psi_2 F_x\|_{{\tilde L}^4_{t,\lambda}}
 \|F\|^{k-1}_{A_{\lambda}} \label{mi2b}
 \end{equation}
 $ \bullet $ {\it Estimate on $ \Lambda_2$. }
 We rewrite $G_2 $ as
\begin{eqnarray*}
G_2 & = &  Q_{\{2^{10}k\}}\Bigl( \sum_{i_2\ge 8+\alpha(k)}
\Delta_{i_2}({\tilde v})
 \sum_{1\le i_3< i_2-6-\alpha(k)}\hspace*{-5mm} \Delta_{i_3} (F)\hspace*{-5mm}
 \sum_{0\le
 i_4,..,i_{k+2}\le i_3}\hspace*{-5mm}n(i_3,..,i_{k+2})
 \prod_{j=4}^{k+2} \Delta_{i_j} (F)\Bigr)  \\
  & & + Q_{\{2^{10}k\}}\Bigl( \sum_{i_2\ge 8+\alpha(k)}
\Delta_{i_2}({\tilde v}) (\Delta_{0} (F))^k \Bigr)\\
 & =&  \sum_{p\ge 1}\Bigl[  Q_{\{2^{10}k\}}
 \Bigl( \sum_{i_2\ge 8+\alpha(k)\atop i_2> p+6+\alpha(k)}
\hspace*{-5mm}\Delta_{i_2}({\tilde v})
  \Delta_{p} (F) \hspace*{-5mm}\sum_{0\le
 i_4,..,i_{k+2}\le p} \hspace*{-5mm}n(i_3,..,i_{k+2})\prod_{j=4}^{k+2} \Delta_{i_j} (F)\Bigr)\Bigr]  \\
  & & + \, Q_{\{2^{10}k\}}\Bigl( \sum_{i_2\ge 8+\alpha(k)}
\Delta_{i_2}({\tilde v}) (\Delta_{0} (F))^k \Bigr)\\
& = & \sum_{p\ge 1} H_p + L \quad .
\end{eqnarray*}
it is thus clear that
\begin{eqnarray*}
 \Lambda_2 & \lesssim & \sum_{p\ge1}
 \Bigl| \int_{B_2}  \xi\,  |\hat{h}(\tau,\xi)|
  \xi_1^{-1} |\hat{w}(\tau_1,\xi_1|
  |\xi-\xi_1||\widehat{H_p}(\tau-\tau_1,\xi-\xi_1)|\\
 & & +
  \int_{B_2}  \xi\,  |\hat{h}(\tau,\xi)|
  \xi_1^{-1} |\hat{w}(\tau_1,\xi_1| |\xi-\xi_1|
  |\hat{L}(\tau-\tau_1,\xi-\xi_1)| \\
  & = & \Lambda_{21}+ \Lambda_{22} \quad .
\end{eqnarray*}
 We  rewrite  $
 \Lambda_{21} $ as the sum of two terms :
\begin{eqnarray*}
 \Lambda_{21} & = &\sum_{p\ge1}\int_{B_2} \chi_{\{|\xi|\le
2^{p+6+\alpha(k)}\}} \xi\, |\hat{h}(\tau,\xi)|
  \xi_1^{-1} |\hat{w}(\tau_1,\xi_1)| |\xi-\xi_1|
|\widehat{H_p}(\tau-\tau_1,\xi-\xi_1)| \\
& & +\sum_{p\ge1}\int_{B_2} \chi_{\{|\xi|> 2^{p+6+\alpha(k)}\}}
\xi\, |\hat{h}(\tau,\xi)|
  \xi_1^{-1} |\hat{w}(\tau_1,\xi_1)| |\xi-\xi_1|
|\widehat{H_p}(\tau-\tau_1,\xi-\xi_1)| \\
& =&  \Lambda_{21}^1+ \Lambda_{21}^2 \quad .
\end{eqnarray*}
Let us explain the idea of this dichotomy. In the domain of
integration of $ \Lambda_{21}^1 $, the frequency $ \xi $ of $
\hat{h} $ is controlled by the maximum of the $|\xi_i|$,
$i=3,..,k+1 $,
 and thus we can in some sens exchange the derivative on $ h$ with a
 derivative on $ F $. On the other hand, in the domain of
 integration of $\Lambda_{21}^2 $, $|\xi| $ and $|\xi_2| $ are
 very large with respect to $ |\xi_3|, .., |\xi_{k+1}|$ and then
 we  have a
 good non-resonant relation (similar to the non-resonant relation
  used in \cite{M1} to prove the  bilinear estimate
  \re{LL1}) that enables to regain one derivative.

 $ \bullet $ {\it Estimate on $
\Lambda_{21}^1 $}.
  Using a Littlewood-Paley decomposition
 of $ h $, we get thanks to \re{C3} and Cauchy-Schwarz inequality in $ p$
 \begin{eqnarray*}
 \Lambda_{21}^1 & \lesssim &
\sum_{p\ge1}\sum_{q=-7-\alpha(k)}^{p-8-\alpha(k)}
  \int_{B_2}  2^{p-q}  |\widehat{\Delta_{p-q} h}(\tau,\xi)|
  |\hat{w}(\tau_1,\xi_1)|
|\widehat{H_p}(\tau-\tau_1,\xi-\xi_1)|\quad \\
 & \lesssim &
 \sum_{q=-7-\alpha(k)}^\infty  2^{-q} \sum_{p\ge q+8+\alpha(k)}
  \int_{B_2}   |\widehat{\Delta_{p-q} h}(\tau,\xi)|
  |\hat{w}(\tau_1,\xi_1)|
2^p|\widehat{H_p}(\tau-\tau_1,\xi-\xi_1)|  \\
 &  \lesssim  &
\sup_{q\ge -7-\alpha(k)}  k \sum_{p\ge q+8+\alpha(k)}
  \int_{B_2}    |\widehat{\Delta_{p-q} h}(\tau,\xi)|
  |\hat{w}(\tau_1,\xi_1)|
2^p|\widehat{H_p}(\tau-\tau_1,\xi-\xi_1)| \\
&   \lesssim & k \| {\cal F}^{-1} (\langle \xi\rangle^{s}|\hat{w}|
)\|_{L^4_{t,\lambda}} \Bigl(\sum_{p\ge 1} \|{\cal F}^{-1} (\langle
\xi\rangle^{-s} |\widehat{\Delta_p h}| )\|_{L^4_{t,\lambda}}^2
\Bigr)^{1/2} \Bigl(\sum_{p\ge 1} 2^{2p}
\|H_p\|_{L^2_{t,\lambda}}^2\Bigr)^{1/2}
\end{eqnarray*}
Note that $ {\tilde L}^4_{t,\lambda}\hookrightarrow X^{3/8,0}  $
since by \re{l5}, for any function $ z\in X^{3/8,0}_\lambda $,
\begin{equation}
\Bigl(\sum_{p\ge 1} \|{\cal F}^{-1} (|\widehat{\Delta_p z}|)
\|_{L^4_{t,\lambda}}^2 \Bigr)^{1/2} \lesssim \Bigl(\sum_{p\ge 1}
\|\Delta_p z\|_{X^{3/8,0}_\lambda}^2 \Bigr)^{1/2} \lesssim
\|z\|_{X^{3/8,0}_\lambda} \; . \label{l4}
\end{equation}
Moreover since, according to the frequency localization of the
functions,
$$
\Delta_q H_p=   Q_{2^{10}k}\Bigl( \sum_{i_2\ge 8+\alpha(k)\atop
q-1\le i_2\le q+1} \hspace*{-5mm}\Delta_{i_2}({\tilde v})
  \Delta_{p} (F)\hspace*{-5mm} \sum_{0\le
 i_4,..,i_{k+2}\le p} \hspace*{-5mm}n(p,i_4,..,i_{k+2})
\prod_{j=4}^{k+2} \Delta_{i_j} (F)\Bigr)
  $$
  we infer from \re{estA} and \re{l4} that
  \begin{eqnarray*}
   \|H_p\|_{L^2_{t,\lambda}}^2 &\sim & \sum_{q\ge p+9+\alpha(k)}
   \|\Delta_{q} H_p\|_{L^2_{t,\lambda}}^2 \\
 & \lesssim &  k^2 \|F\|_{A_{\lambda}}^{2(k-1)} \|\psi_2 \Delta_p F\|_{L^4_{t,\lambda}}^2
  \sum_{q\ge 1} \|\Delta_q {\tilde v}\|_{L^4_{t,\lambda}}^2 \\
 &  \lesssim & k^2 \|F\|_{A_{\lambda}}^{2(k-1)} \|v\|_{X^{1/2,0}}^2
  \|\psi_2 \Delta_p F \|_{L^4_{t,\lambda}}^2 \quad .
\end{eqnarray*}
 Therefore, we deduce that
 \begin{eqnarray}
\Lambda_{21}^1 & \lesssim & k
 \|h\|_{X^{3/8,-s}_\lambda} \|w\|_{X^{1/2,s}_\lambda}  \|F\|_{A_{\lambda}}^{k-1} \|v\|_{X^{1/2,0}_\lambda}
 \Bigl( \sum_{p\ge 1} 2^{2p} \| \psi_2 \Delta_p F \|_{{\tilde
 L}^4_{t,\lambda}}^2 \Bigr)^{1/2}\nonumber  \\
& \lesssim & k \,
 \|h\|_{X^{3/8,-s}_\lambda} \|w\|_{X^{1/2,s}_\lambda}  \|F\|_{A_{\lambda}}^{k-1} \|v\|_{X^{1/2,0}_\lambda}
 \|\psi_2 F_x \|_{{\tilde L}^4_{t,\lambda}} \quad .\label{mi3}
\end{eqnarray}
$ \bullet $ {\it Estimate on $ \Lambda_{21}^2 $ and $ \Lambda_{22}
$.}  Since clearly ,  $ \sum_{p\ge 0} |\widehat{\Delta_p (f)
}(\tau,\xi)|\le 2  |\hat{f}(\tau,\xi)| $ for any $ f\in
L^2_{t,\lambda} $, we infer that
 \begin{eqnarray*}
\Lambda_{21}^2 & \lesssim &  k \sum_{p\ge 1} \sum_{i_2\ge
p+7+\alpha(k)}
  \int_{B_3}  \xi |\widehat{ h}(\tau,\xi)|\xi_1^{-1}
  |\hat{w}(\tau_1,\xi_1)| |\xi-\xi_1||\widehat{\Delta_{i_2}({\tilde
  v})}(\tau_2,\xi_2)|
  |\widehat{\Delta_{p} (F)}(\tau_3,\xi_3)| \\
   & &\hspace*{20mm}  \sum_{
 i_4,..,i_{k+2}\ge 0} \prod_{j=4}^{k+2} |\widehat{\Delta_{i_j}
 (F)}(\tau_j,\xi_j)|  \\
& \lesssim & k \int_{B_3}\xi\, |\hat{h}(\tau,\xi)|
  \xi_1^{-1} |\hat{w}(\tau_1,\xi_1)|
  |\xi-\xi_1||\hat{v}(\tau_2,\xi_2)| \prod_{i=3}^{k+2}
  |\hat{F}(\tau_i,\xi_i)| \\
  & =& {\tilde J}_{k_{/B_3}}
\end{eqnarray*}
where
$$
 B_3=\{(\tau,\tau_1,..,\tau_{k+1}, \xi,\xi_1,..,\xi_{k+1})\in B_1,\,
 \, \xi_2\le - 2^{10}k,\, \min(|\xi|,|\xi_2|)> 10 k \, \max_{i=3,..,k+2} |\xi_i| \}\quad .
$$
In the same way, it is easy to check that
  $ \Lambda_{22} \lesssim {\tilde  J}_{k_{/B_3}} $. We set  $
\sigma=\sigma(\tau,\xi)=\tau-\xi|\xi| $ and
$\sigma_i=\sigma(\tau_i,\xi_i) $, $ i=1,..,k+2 $. Noticing that
on $B_3 $ the sign
 of $ \xi,\,\xi_1 $ and $\xi_2 $ are known,   we  get  the following
algebraic relation :
\begin{eqnarray}
\sigma-\sum_{i=1}^{k+2} \sigma_i & = &  (\sum_{i=1}^{k+2} \xi_i)^2-\xi_1^2+\xi_2^2-\sum_{i=3}^{k+2}\xi_i|\xi_i|\nonumber \\
& = & 2 \xi_2  \sum_{i=1}^{k+2} \xi_i+2\xi_1\sum_{i=3}^{k+2} \xi_i -
 \sum_{i=3}^{k+2}\xi_i|\xi_i|+\Bigl(\sum_{i=3}^{k+2} \xi_i\Bigr)^2\nonumber\\
& = & 2 \xi_2  \xi +2 \xi_1\sum_{i=3}^{k+2} \xi_i -
 \sum_{i=3}^{k+2}\xi_i|\xi_i|+\Bigl(\sum_{i=3}^{k+2} \xi_i\Bigr)^2\quad .\label{mi5}
\end{eqnarray}
Note that  on $B_3 $, we have  $(\sum_{i=3}^{k+2} \xi_i)^2 \le
10^{-2} |\xi_2 \xi | $, $ \sum_{i=3}^{k+2}\xi_i|\xi_i| \le 10^{-2}
|\xi_2 \xi | $  and
\begin{equation}
\frac{1}{2} |\xi_2|\le  |\xi-\xi_1| \le
|\xi_2|+\sum_{i=3}^{k+2}|\xi_i|\le
 2 |\xi_2| \label{mi6} \quad .
\end{equation}
Hence, $ \xi_1 \le 2\max (|\xi|,|\xi-\xi_1|) \le 4\max
(|\xi|,|\xi_2|) \, $ and  $| \xi_1  \sum_{i=3}^{k+2}\xi_i| \le
2|\xi_2 \xi |/5 $.  We thus
  deduce from \re{mi5} that the following non-resonant relation
  holds
\begin{equation}
\max_{i=1,..,k+2}(|\sigma|,|\sigma_i|) \gtrsim |\xi \xi_2 |/k  \quad .\label{mi7}
\end{equation}
It remains to divide $ B_3 $ in subregions according to the indice
where the maximum
  is reached in \re{mi7}. Thanks to \re{C3},
$$
{\tilde J}_{k_{/B_3}}  \lesssim  \int_{B_2} |\xi \xi_2|^{1/2}
|\hat{h}(\tau,\xi) |
  | \hat{w}(\tau_1,\xi_1)|
 |\hat{v}(\tau_2,\xi_2)|  \prod_{i=3}^{k+2}|\hat{F}(\tau_i,\xi_i)|\quad .
$$
$\bullet \, |\sigma| $ dominant. By \re{C3} and \re{mi7},
Plancherel, Holder inequality and \re{l5}, we infer that
\begin{eqnarray}
{\tilde J}_{k_{/B_3}}  & \lesssim & k \| {\cal F}^{-1} (\langle
\sigma\rangle^{1/2}\langle \xi\rangle^{-s}|\hat{h}|
)\|_{L^2_{t,\lambda}} \| {\cal F}^{-1} (\langle
\xi\rangle^{s}|\hat{w}| )\|_{L^4_{t,\lambda}} \| {\cal F}^{-1}
(|\hat{v}| )\|_{L^4_{t,\lambda}}
\| {\cal F}^{-1} (|\hat{F}| )\|_{L^\infty_{t,\lambda}}^{k} \nonumber\\
& \lesssim &  k \|h\|_{X^{1/2,-s}_\lambda}\|w\|_{X^{1/2,s}_\lambda}\|v\|_{X^{1/2,0}_\lambda}
 \|F\|_{A_\lambda}^k \label{mi8}\quad .
\end{eqnarray}
 $\bullet \, |\sigma_1| $ or $|\sigma_2|  $ dominant. It
is easy to see that in the same way
 \begin{eqnarray}
{\tilde J}_{k_{/B_3}} & \lesssim &
  k \|h\|_{X^{3/8,-s}_\lambda}\|w\|_{X^{1/2,s}_\lambda}\|v\|_{X^{1/2,0}_\lambda}
 \|F\|_{A_\lambda}^k \label{mi9}\quad .
\end{eqnarray}
$\bullet \, |\sigma_i|, \, i\ge 3 $,
  dominant. By Plancherel, Holder inequality and \re{l5}, we infer that
\arraycolsep1pt
\begin{eqnarray}
{\tilde J}_{k_{/B_3}}  & \lesssim & k^{1/2} \| {\cal F}^{-1}
(\langle \xi\rangle^{-s}|\hat{h}| )\|_{L^4_{t,\lambda}}
\| {\cal F}^{-1} (\langle \xi\rangle^{s}|\hat{w}| )\|_{L^4_{t,\lambda}} \nonumber \\
& & \quad\| {\cal F}^{-1} (|\hat{v}| )\|_{L^4_{t,\lambda}}
\| {\cal F}^{-1}(| \sigma|^{1/2} \chi_{\{|\sigma|\gtrsim 1\}}|\hat{F}| )\|_{L^4_{t,\lambda}}
\| {\cal F}^{-1}( |\hat{F}| )\|_{L^\infty_{t,\lambda}}^{k-1}
 \nonumber \\
& \lesssim &  k  \|h\|_{X^{3/8,-s}_\lambda}\|w\|_{X^{1/2,s}_\lambda}
\|v\|_{X^{1/2,0}_\lambda}
 (\|P_3 F\|_{{\dot X}^{1,0}_\lambda}+
 \|P_{>3} F \|_{X^{7/8,0}_\lambda})\|F\|_{A_\lambda}^{k-1}\label{mi10} \quad ,
\end{eqnarray}
\arraycolsep5pt
 where  we use that
 $ |\sigma_i|\gtrsim |\xi\xi_2|/k \gtrsim 1 $ on $ B_3 $  to get an homogeneous
  Bourgain type norm on $ P_3 F $. It has some importance when
    using dilations argument since the $ L^2 $-norm
   of $ P_3 F $ is surcritical and thus behaves badly for such arguments.
    Since clearly, $\|P_{>3} F \|_{X^{7/8,0}_\lambda}\lesssim
     \|P_{>3} F_x\|_{X^{7/8,-1}_\lambda} $,
 gathering \re{mi8bbb}, \re{mi8bb}, \re{mi2b}, \re{mi3}, \re{mi8}, \re{mi9} and \re{mi10}, \re{mi2}
follows.
\begin{lem} For any $ s\ge 0 $ it holds
\begin{eqnarray}
\hspace*{-10mm}\Bigl\|  \partial_x P_+ \Bigl( W\partial_x
P_-(e^{-iF/2} P_- v )\Bigr)\Bigr\|_{Z^{-1,s}_\lambda}
 & \lesssim &   \|w\|_{X^{1/2,s}_\lambda}
 \|v\|_{X^{1/2,0}_\lambda}\, e^{C\|F\|_{A_\lambda}}\nonumber \\
 & & \hspace*{-70mm}\Bigl(1+ \|P_3 F\|_{{\dot X}^{1,0}_\lambda} +\|P_{>3}
F_x \|_{X^{7/8,-1}_\lambda} + \|F\|_{A_\lambda} + \|\psi_2 F_x
\|_{\tilde{L^4}_{t,\lambda}} \Bigr)\; .\label{mi11}
\end{eqnarray}
\label{nonlin3}
\end{lem}
{\it Proof.} The proof of this lemma essentially follows  the one
of Lemma \ref{nonlin2} and thus will be only sketched. We estimate
$$
  \Bigl\|  \partial_x P_+ \Bigl( W\partial_x  P_- v )\Bigr)\Bigr\|_{Z^{-1,s}_\lambda}
 +\sum_{k\ge 1}\ \frac{1}{k!}
 \Bigl\|  \partial_x P_+ \Bigl( W\partial_x P_-(F^k P_- v )\Bigr)\Bigr\|_{Z^{-1,s}_\lambda} \; .
$$
Again, the first term of the above inequality is estimated in
(\cite{M1}, Lemma 3.4) by
\begin{equation}
\Bigl\|  \partial_x P_+ \Bigl( W\partial_x  P_- v
)\Bigr)\Bigr\|_{Z^{-1,s}_\lambda} \lesssim
\|w\|_{X^{1/2,s}_\lambda} \|v\|_{X^{1/2,0}_\lambda}
\end{equation}
To estimate the second term we
 first note that   by Cauchy-Schwarz in $ \tau $,
$$
\Bigl\|\frac{\langle \xi \rangle^{s}}{ \langle \sigma \rangle} {\cal F}
\Bigl[\partial_x {P}_+ \Bigl( W\partial_x P_-(F^k P_- v )\Bigr)\Bigr]  \Bigr\|_{L^2_\xi L^1_\tau} \lesssim
\Bigl\|\partial_x {P}_+   \Bigl( W\partial_x P_-(F^k P_- v )\Bigr)\Bigr\|_{X^{-1/2+\varepsilon,s}_\lambda}, \; \varepsilon>0 .
$$
On account of \re{mi8bbb}, \re{mi8bb}, \re{mi2b}, \re{mi3}, \re{mi9} and \re{mi10}, this last term is controlled
  by the right-hand side of \re{mi11} except in the region $B_3 $
  with  $ |\sigma| $ dominant. \\
 Moreover, in the region $ \{\xi_1\le 1 \} $, using \re{C3} and then \re{l5} we infer that
\begin{eqnarray*}
\Bigl\|\partial_x {P}_+\Bigl( W\partial_x P_-(F^k P_- v )\Bigr)
\Bigr\|_{X^{-1/2+\varepsilon,s}_\lambda}
   & \lesssim  & \|{\cal F}_{t,x}^{-1}(|w|)\|_{L^4_{t,\lambda }}
    \|{\cal F}_{t,x}^{-1}(|v|)\|_{L^4_{t,\lambda }}
 \|{\cal F}_{t,x}^{-1}(|F|)\|_{L^\infty_{t,\lambda }}^k \nonumber \\
    & \lesssim & \|w\|_{X^{1/2,0}_\lambda} \, \|v\|_{X^{1/2,0}_\lambda}
\|F \|_{A_\lambda }^k
  \quad .
\end{eqnarray*}
It thus remains to treat the region $B_3 $ with $\xi_1\ge 1 $ and
$ |\sigma| $ dominant. To handle with this region we proceed as in
\cite{CKSTT2}. The proof is very similar to the one
of Lemma 3.4 in \cite{M1}. \\
By \re{mi5} in this region we have
\begin{equation}
\langle\sigma \rangle \gtrsim \langle \xi \xi_2 \rangle /
k\gtrsim 1 \label{mi12}\quad .
\end{equation}
Therefore \re{mi11} will be proven if we show the following inequality:
\begin{equation}
J_k \lesssim  k \| \tilde{w}\|_{L^2_{t,\lambda}} \|
\tilde{v}\|_{L^2_{t,\lambda}} \| \hat{F}\|_{L^1_{\tau,\xi} }^k
\label{mi13}
\end{equation}
with
\begin{equation}
J_k  =\Bigl\|\int_{C(\tau,\xi)} \frac{\langle \xi \rangle^{s} \xi
\langle \xi_1 \rangle^{-s}\xi_1^{-1}
|\widehat{\tilde{w}}(\tau_1,\xi_1)||\xi-\xi_1|
|\widehat{\tilde{v}}(\tau_2,\xi_2)|
\displaystyle \prod_{i=3}^{k+2}|\hat{F}(\tau_i,\xi_i)|}{\langle \sigma \rangle
\langle \sigma_1 \rangle^{1/2}
 \langle \sigma \rangle^{1/2}}\Bigr\|_{L^2_\xi L^1_\tau} \label{mi14}
\end{equation}
and
\begin{eqnarray}
C(\tau,\xi)& = & \Bigl\{(\tau_1,..,\tau_{k+1}, \xi_1,..,\xi_{k+1})\in \R^{k+1}\times(\lambda^{-1}\Z)^{k+1}, \nonumber \\
& &(\tau,\tau_1,..,\tau_{k+1},\xi, \xi_1,..,\xi_{k+1}) \in B_3,\,
|\xi_1|>1, \max_{i=1,..,k+2} (\langle \sigma_i \rangle)\le
\langle \sigma \rangle
 \Bigr\} \quad. \label{mi15}
\end{eqnarray}
$ \bullet $  The subregion $  \max(|\sigma_1|,|\sigma_2|) \ge (\xi |\xi_2|)^{\frac{1}{16}} $.
 We will assume that $ \max(|\sigma_1|,|\sigma_2|)=|\sigma_1| $ since the other case can be treated in exactly the same way.
By \re{mi12} and \re{mi6}, recalling that on the domain of integration $ \xi_1 \ge \max(\xi, |\xi-\xi_1|) $,  we infer that
$$
J_k  \lesssim  k \Bigl\|  \int_{C_1(\tau,\xi)}
\frac{|\widehat{\tilde{w}}(\tau_1,\xi_1)|
|\widehat{\tilde{v}}(\tau_2,\xi_2)|
\displaystyle\prod_{i=3}^{k+2}|\hat{F}(\tau_i,\xi_i)|}{\langle \sigma
\rangle^{1/2+\frac{1}{128}} \langle\sigma_1 \rangle^{3/8}
 \langle \sigma \rangle^{1/2}}\Bigr\|_{L^2_\xi L^1_\tau}
$$
where
$$ C_1(\tau,\xi)=\{(\tau_1,\xi_1)\in C(\tau,\xi), \,  |\sigma_1| \ge (\xi |\xi_2|)^{\frac{1}{16}} \} \quad .$$
Applying Cauchy-Schwarz in $ \tau $ we obtain thanks to \re{l5},
\begin{eqnarray*}
J_k   & \lesssim & k \Bigl\| \int_{C_1(\tau,\xi)}
\frac{|\widehat{\tilde{w}}(\tau_1,\xi_1)|
|\widehat{\tilde{v}}(\tau_2,\xi_2)|
\displaystyle\prod_{i=3}^{k+2}|\hat{F}(\tau_i,\xi_i)|}{ \langle\sigma_1
\rangle^{3/8}
 \langle \sigma \rangle^{1/2}}\Bigr\|_{L^2_{\xi,\tau}}\nonumber \\
    &  \lesssim &  k \Bigl\|{\cal F}^{-1}\Bigl( \frac{|\widehat{\tilde{w}}|}
{\langle\sigma \rangle^{3/8}}\Bigr) \Bigr\|_{L^4_{t,\lambda}}\,
\Bigl\|{\cal F}^{-1}\Bigl( \frac{|\widehat{\tilde{v}}|}{\langle\sigma \rangle^{1/2}}\Bigr)\Bigr\|_{L^4_{t,\lambda}}
 \|{\cal F}^{-1}(|\hat{F}|)\|_{L^\infty_{t,\lambda}}^k
\nonumber \\
     & \lesssim &  k \|\tilde{w}\|_{L^2_{t,\lambda}}\,
      \|\tilde{v}\|_{L^2_{t,\lambda}} \|\hat{F}\|_{L^1_{\tau,\xi}}^k \quad .
\end{eqnarray*}
$\bullet $ The subregion $  \max(|\sigma_1|,|\sigma_2|)
 \le (\xi |\xi_2|)^{\frac{1}{16}} $. Changing the $ \tau, \tau_1,.., \tau_{k+1} $
  integrals in $ \tau_1,..,\tau_{k+2} $ integrals in \re{mi14} and using \re{mi6}
   and \re{mi12},
   we infer that
\begin{eqnarray*}
J_k  & \lesssim & k \Bigl\|
\chi_{\{\xi\ge 1\}} \int_{D(\xi)} \xi_1^{-1}\int_{\tau_1=-\xi_1^2+O(|\xi \, \xi_2|^{1/16})} \frac{|\widehat{\tilde{w}}(\tau_1,\xi_1)|}{\langle \tau_1+|\xi_1|\xi_1\rangle^{1/2}} \\
 & & \quad\quad \int_{\tau_2=\xi_2^2+O(|\xi \, \xi_2|^{1/16})}
\frac{|\widehat{\tilde{v}}(\tau_2,\xi_2)|}{\langle \tau_2+|\xi_2|
\xi_2\rangle^{1/2}} \int_{\tau_3,..,\tau_{k+2}}
\prod_{i=3}^{k+2}|\hat{F}(\tau_i,\xi_i)|\Bigr\|_{L^2_\xi}
\end{eqnarray*}
with
$$
 D(\xi)=\{(\xi_1,..,\xi_{k+1})\in (\lambda^{-1}\Z)^{k+1}, \xi_1\ge 1, \,
\xi-\xi_1\le -1/\lambda \,  \} \; .
$$
 Applying Cauchy-Schwarz inequality in $
\tau_1 $ and $ \tau_2 $ and recalling that $ \xi_1\ge 1 $ we get
$$
J_k \lesssim k \Bigl\|\chi_{\{\xi\ge 1\}} \int_{D(\xi)} \langle
\xi_1\rangle^{-1} (\xi |\xi_2|)^\frac{1}{8} K_1(\xi_1) K_2(\xi_2)
 \prod_{i=3}^{k+2}K(\xi_i)\Bigr\|_{L^2_\xi}
$$
where
$$
K_1(\xi)=\Bigl( \int_\tau \frac{|\widehat{\tilde{w}}(\tau,\xi)|^2}{\langle \tau+|\xi|\xi\rangle} \Bigr)^{1/2}, \quad
K_2(\xi)=\Bigl( \int_\tau \frac{|\widehat{\tilde{v}}
(\tau,\xi)|^2}{\langle \tau+|\xi| \xi \rangle} \Bigr)^{1/2} \mbox{ and }
 K(\xi)=\int_\tau |\hat{F}(\tau,\xi)| \; .
$$
Therefore, by using \re{mi6} and \re{C3}, H\"older and then
Cauchy-Schwarz inequalities,
\begin{eqnarray}
J_k  & \lesssim & k \Bigl\|  \langle\xi\rangle^{-\frac{3}{4}}
\int_{(\xi_3,..,\xi_{k+2})\in (\lambda^{-1} \Z)^k}\prod_{i=3}^{k+2} K (\xi_i)  \int_{\xi_1\in \lambda^{-1}
\Z} K_1(\xi_1) K_2(\xi_2) \Bigr\|_{L^2_\xi} \nonumber \\
 & \lesssim & k \Bigl\|\int_{(\xi_3,..,\xi_{k+2})\in (\lambda^{-1} \Z)^k}\prod_{i=3}^{k+2} K (\xi_i)
 \int_{\xi_1\in\lambda^{-1}\Z} K_1(\xi_1) K_2(\xi_2)
\Bigr\|_{L^\infty_\xi} \nonumber \\
  & \lesssim & k \int_{(\xi_3,..,\xi_{k+2})\in (\lambda^{-1} \Z)^k}\prod_{i=3}^{k+2} K (\xi_i)
   \Bigl(\int_{\xi \in \lambda^{-1}\Z} K_1(\xi)^2 \Bigr)^{1/2}\, \Bigl(\int_{\xi_\in
   \lambda^{-1}\Z} K_2(\xi)^2 \Bigr)^{1/2} \nonumber \\
   &  \lesssim & k \|\tilde{w}\|_{X^{-1/2,0}_{\lambda}} \, \|\tilde{v}\|_{X^{-1/2,0}_{\lambda}} \|\hat{F}\|_{L^1_{\tau,\xi}}^k\nonumber\\
   & \lesssim & \, k  \|\tilde{w}\|_{L^2_{t,\lambda}} \, \|\tilde{v}\|_{L^2_{t,\lambda}}
 \|\hat{F}\|_{L^1_{\tau,\xi}}^k \quad ,\label{C12}\vspace*{3mm}
\end{eqnarray}

\subsection{ End of the proof of \re{nonlinear1}.}\label{Sec42}
 It remains to treat the third term  of the right-hand side of \re{eq2w}. Observe that
 by Cauchy-Schwarz inequality in $ \tau $, Sobolev inequalities in time and Minkowski inequality,
 $$
\| P_0(u^2) w \|_{Z_\lambda^{-1,s}} + \| P_0(u^2) w \|_{X^{-1/2,s}_\lambda}   \lesssim
 \| P_0(u^2) w \|_{X^{-1/2+\varepsilon',s}_\lambda}
 \lesssim \| P_0(u^2) w \|_{L^{1+\varepsilon}_t H^s_\lambda} \quad,
 $$
 for some $ 0<\varepsilon,\varepsilon'<\!\!<1 $.
Assuming that  $ w$ is supported in time in $ [-2,2] $, by H\"older inequality in time and \re{l5} we get
 $$
\| P_0(u^2) w \|_{L^{1+\varepsilon}_t H^s_\lambda}  \lesssim   \|J^s_x w \|_{L^4_{t,\lambda}} \|\psi_2^2 P_0(u^2) \|_{L^2_t L^4_\lambda}
   \lesssim   \|w\|_{X^{1/2,s}_\lambda} \|\psi_2^2 P_0(u^2) \|_{L^2_{t,\lambda}} \quad ,
$$
where we used that  $ \|1\|_{L^4_\lambda} \le  \|1\|_{L^2_\lambda} $ since $ \lambda\ge 1 $.
Hence, the following estimate holds:
\begin{equation}
\| P_0(u^2) w \|_{Z_\lambda^{-1,s}} + \| P_0(u^2) w \|_{X^{-1/2,s}_\lambda}   \lesssim
   \|w\|_{X^{1/2,s}_\lambda} \|\psi_2 u \|^2_{L^4_{t,\lambda}} \quad .\label{termsup}
\end{equation}
Therefore,
 combining Lemmas \ref{line2}-\ref{line3},
\ref{nonlin1}-\ref{nonlin2}-\ref{nonlin3} and \re{termsup},
  we infer that for $  s\ge 0 $, the extension
 $ w^* $ of $ w$ defined by \re{wetoile} satisfies
 \begin{eqnarray*}
 \|w^*\|_{Y^s_\lambda} & \lesssim & \|w(0)\|_{H^s_\lambda} + \|\tilde W_x\|_{X^{1/2,s}_\lambda}
  \, e^{C\, \|{\tilde F}\|_{A_\lambda} }\Bigl[ \|\tilde{u}\|^2_{L^4_{t,\lambda}}
  + \|\tilde{\tilde w}\|_{X^{1/2,0}_\lambda} \nonumber \\
 & &   \hspace*{-10mm} \Bigl( \|P_3 \tilde{F}\|_{{\dot X}^{1,0}_\lambda}+
\|P_{>3} \tilde{F}_x\|_{X^{7/8,-1}_\lambda}+ \|\tilde{F}
\|_{A_\lambda} +
\|\psi_2 \tilde{F}_x \|_{\tilde{L^4}_{t,\lambda}}\Bigr) \Bigr]  \nonumber \\
  & & \hspace*{-25mm} \lesssim
\|w(0)\|_{H^s_\lambda} + \|\ w \|_{X^{1/2,s}_{1,\lambda}} \Bigl(
 \|u\|^2_{N_{1,\lambda}}+\|\ w \|_{X^{1/2,0}_{1,\lambda}}\Bigr)
  e^{C\, ( \|\widehat{\partial_x^{-1} u_0} \|_{L^1_\xi}
  +\|u\|_{N_{1,\lambda}}+\|u\|^2_{N_{1,\lambda}})  }
 \end{eqnarray*}
where in the last step we used  Lemma \ref{estimbasfreq} below
 to estimate $ \|P_3 \tilde{F}\|_{{\dot X}^{1,0}_\lambda} +
    \|P_3 \tilde{F} \|_{A_\lambda} $ and
that, by Cauchy-Schwarz
 in
 $ \xi $,
\begin{equation}
  \| P_{>3} \tilde{F}\|_{A_\lambda}
  \lesssim \|\widehat{{\tilde F}_x}\|_{L^2_\xi L^1_\tau} \quad .
  \label{itit}
  \end{equation}
  \begin{lem}
\label{estimbasfreq}
Let $ {\tilde {\tilde u}} \in N_{1,\lambda} $  and let $ P_3 \tilde{u} $ be defined as in
\re{eqbasfreq}. Then $ P_3 {\tilde F}= P_3 \partial_x^{-1} {\tilde u} $ satisfies :
\begin{equation}
\|P_3 {\tilde F} \|_{A_\lambda} \lesssim \| \widehat{
\partial_x^{-1}  u_0} \|_{L^1_\xi} + \|{\tilde
{\tilde u}}\|^2_{N_\lambda}\label{lala1} \quad ,
\end{equation}
\begin{equation}
\|P_3 {\tilde F}\|_{\dot{X}^{1,0}_\lambda}
  \lesssim  \|u_0
\|_{L^2_\lambda}+\|{\tilde{ \tilde u}}\|^2_{N_\lambda}
\label{lala2}
\end{equation}
and
\begin{equation}
\|P_3 {\tilde F}_x \|_{N_\lambda}
  \lesssim \|u_0
\|_{L^2_\lambda} + \|{\tilde{ \tilde u}}\|^2_{N_\lambda} \label{lala3}\quad .
\end{equation}
Moreover, $ \forall\;  0<\alpha<3 $,
\begin{equation}
\|\psi_2  P_3 Q_{\alpha} {\tilde F} \|_{X^{7/8,0}_\lambda} \lesssim \frac{1}{\alpha} \|u_0\|_{L^2_\lambda}
 +\|{\tilde {\tilde u}}\|_{N_\lambda}^2 \quad . \label{lala4}
\end{equation}
\end{lem}
We postponed the proof of this lemma to the end of this section. \vspace*{2mm}\\
On the other hand, obviously,
$$
\|P_{>1} w^* \|_{X^{1,-1}_\lambda} \lesssim
 \|P_{>1} \partial_x^{-1} w^* \|_{X^{1,0}_\lambda}
 $$
and from \re{wetoile} we deduce that $ w^*=\psi w^{**} $ where $
w^{**} $ satisfies \re{eq2w} with $ W $, $ w$ and $ F $
respectively replaced by $ \psi {\tilde W}$, $ \psi
{\tilde{\tilde w}} $ and $ {\tilde F }$ in the right-hand side member. Therefore using Lemma
\ref{non2} and expanding the exponential function we infer that
\arraycolsep1pt
\begin{eqnarray*}
\|P_{>1} w^* \|_{X^{1,-1}_\lambda} & \lesssim &
\|w^*\|_{L^\infty_t L^2_\lambda} + \|\partial_t (P_{>1}
\partial_x^{-1} w^{**})+ {\cal H} \partial_x^2  (P_{>1}
\partial_x^{-1}
w^{**})\|_{L^2_{t,\lambda}} \nonumber \\
 & \lesssim & \|w^*\|_{Y^s_\lambda}+\|\tilde{w}\|_{L^4_{t,\lambda}}\Bigl(
\|\tilde{\tilde w}\|_{L^4_{t,\lambda}} + \|\psi_2 \tilde{F}_x
\|_{L^4_{t,\lambda}}+\|\psi_2
\tilde{F}_x\|_{L^4_{t,\lambda}}^2\Bigr)
e^{C\|\tilde{F}\|_{A_\lambda}} \nonumber \\
 && \hspace*{-25mm} \lesssim  \|w(0)\|_{H^s_\lambda} +\|\ w \|_{X^{1/2,s}_{1,\lambda}} \Bigl(
 \|u\|_{N_{1,\lambda}}+\|\ w \|_{X^{1/2,0}_{1,\lambda}}\Bigr)
  e^{C\, ( \|\widehat{\partial_x^{-1} u_0} \|_{L^1}
  +\|u\|_{N_{1,\lambda}}+\|u\|^2_{N_{1,\lambda}})  }
\end{eqnarray*}
\arraycolsep5pt
 Finally, using  Lemma \ref{non1}
 we infer that  for $ 0\le s \le 1/2 $,
\begin{eqnarray}
\|w(0)\|_{H^s_\lambda} & = &
 \|\partial_x P_+ e^{-i\partial_x^{-1} u_0} \|_{H^s_\lambda} = \frac{1}{2}
 \|  P_+ (u_0 e^{-i\partial_x^{-1} u_0})  \|_{H^s_\lambda} \nonumber \\
  & \lesssim & \sum_{k\ge 0} \frac{1}{k!} \|u_0 (\partial_x^{-1} u_0)^k \|_{H^s_\lambda} \nonumber \\
   & \lesssim & \|u_0\|_{H^s_\lambda} (1+\|u_0\|_{L^2_\lambda}) e^{\|\partial_x^{-1}
    u_0 \|_{L^\infty_\lambda}}
\end{eqnarray}
which ends the proof of \re{nonlinear1}.\vspace*{2mm} \\
{\it Proof of Lemma \ref{estimbasfreq}. }
 From \re{eqbasfreq},  $ P_3 \tilde{F} =P_3 \tilde{F}^1+P_3 \tilde{F}^2$
  where  $ P_3 \tilde{F}^1= \psi(t/\lambda^2) V(t) P_3 \partial_x^{-1} u_0 $
 and
   $ P_3 \tilde{F}^2=\psi(t) P_3\tilde{ \tilde{F}^2} $ with
\begin{equation}
P_3\tilde{\tilde{F}^2}(0)=0 \mbox{ and } P_3 \tilde{\tilde{F}^2}_t
+{\cal H}
\partial_x^2 P_3 \tilde{\tilde{F}^2}=P_3\Bigl[(\psi {\tilde
{\tilde u}})^2/2-
 P_0((\psi {\tilde {\tilde u}})^2)/2\Bigr]  \quad . \label{bobo}
\end{equation}
Therefore
$$
\|P_3 {\tilde{F}^2}\|_{\dot{X}^{1,0}_\lambda} \lesssim \|\psi_t
\|_{L^2_t} \|V(t)P_3 \tilde{\tilde{F}^2}\|_{L^\infty_t
L^2_\lambda}+ \|\psi\|_{L^\infty_t} \|P_3 \tilde{\tilde{F}^2}_t
+{\cal H}
\partial_x^2 P_3 \tilde{\tilde{F}^2}\|_{L^2_{t,\lambda}}
$$
 and  \re{bobo}   leads to
\begin{equation}
\| P_3 {\tilde F}^2\|_{\dot{X}^{1,0}_1} \lesssim
  \| \psi {\tilde{\tilde u}}\|_{L^4_{t,\lambda}}^{2}\lesssim \| {\tilde{ \tilde u}}\|^2_{N_\lambda}
  \ \label{otot}
  \end{equation}
On the other hand, by the definition of $ P_3 \tilde{F}^1 $,
\begin{eqnarray}
\| P_3 {\tilde F}^1\|_{\dot{X}^{1,0}_1} & =& \Bigl\|\partial_t
\Bigl( \psi(\cdot/\lambda^2) P_3
\partial_x^{-1} u_0\Bigr) \Bigr\|_{L^2_{t,\lambda}}\nonumber \\
& \lesssim & \lambda^{-2} \|\psi_t(\cdot/\lambda^2) \|_{L^2_t}
\|\partial_x^{-1} u_0\|_{L^2_\lambda} \nonumber\\
& \lesssim & \lambda^{-1}\, \lambda \|u_0\|_{L^2_\lambda}\lesssim
\|u_0\|_{L^2_\lambda}
  \quad . \label{juju}
  \end{eqnarray}
Moreover, from \re{bobo} and Lemma \ref{line2}-\ref{line3} we
deduce that
$$
\|\widehat{P_{3} {\tilde F}}\|_{L^1_{\tau,\xi}} \lesssim
\|\widehat{ P_ 3F_0} \|_{L^1_\xi} +\Bigl\| \chi_{\{|\xi|\le 3\}}
\frac{{\widehat {\psi \tilde {\tilde u}}}\ast{\widehat {\psi
\tilde {\tilde u}}}}
{\langle\sigma\rangle}\Bigr\|_{L^1_{\tau,\xi}}+\Bigl\|\frac{
{\cal F}_t\Bigl(P_0((\psi {\tilde {\tilde u}})^2)\Bigr)}
{\langle\sigma\rangle}\Bigr\|_{L^1_{\tau}}
$$
Applying Cauchy-Schwarz inequality in $ \tau $ and $ \xi $ , it
follows that
\begin{eqnarray*}
\|\widehat{P_{3} {\tilde F}}\|_{L^1_{\tau,\xi}} & \lesssim  & \|
\widehat{P_3 F_0} \|_{L^1_\xi} +\Bigl\| {\widehat {\psi \tilde
{\tilde u}} }\ast {\widehat {\psi \tilde {\tilde u}}}
\Bigr\|_{L^2_{\tau,\xi}}+ \|P_0((\psi {\tilde {\tilde u}})^2)\|_{L^2_t} \\
 & \lesssim & \| \widehat{P_3 F_0}\|_{L^1_\xi}
 +\|\psi {\tilde{\tilde u}}\|_{L^4_{t,\lambda}}^2 \\
& \lesssim & \| \widehat{P_3 F_0}\|_{L^1_\xi}
 + \| {\tilde {\tilde u}}\|_{N_{\lambda}}^2 \quad .
\end{eqnarray*}
To get \re{lala4} we notice that by classical linear estimates in Bourgain spaces (cf. \cite{G}) and \re{eqbasfreq}
 we have
$$
\|\psi_2  P_3 Q_{\alpha} {\tilde F} \|_{X^{7/8,0}_\lambda} \lesssim  \|Q_{\alpha} \partial_x^{-1} u_0\|_{L^2_\lambda}
 +\|(\psi{\tilde {\tilde u}})^2\|_{X^{-1/8,0}_\lambda}
\lesssim \frac{1}{\alpha} \|u_0\|_{L^2_\lambda} +\|\psi {\tilde {\tilde u}}\|_{L^4_{t,\lambda}}^2 \; .
$$
It remains to get the  estimate \re{lala3} on $ \|P_3 {\tilde F}_x
\|_{N_\lambda} $. But this is straightforward  by combining
\re{eqbasfreq}, Lemmas \ref{line2}-\ref{line3} and Sobolev
inequality in time for evaluating the $Z^{0,0}_\lambda $-norm and
by writing $ \|\chi_{[-4,4]}(t) P_3 {\tilde F}_x \|_{{\tilde
L}^4_{t,\lambda}} \lesssim \|\chi_{[-4,4]}(t) P_3 {\tilde F}_x
\|_{L^\infty_t L^2_\lambda} $ and then using the unitarity of $
V(t) $ in $L^2_\lambda $ (see \re{estP1u} below).
\section{Proof of the estimates on $ u $  } \label{Sec5}
In this section we prove estimates \re{nonlinear2} and
\re{nonlinear4} of Proposition \ref{NON}. We will need the
following lemma, proven in the appendix,
 which  enables to  treat the multiplication with
 the gauge function $ e^{-iF/2} $
 in $ {\tilde L}^4_{t,\lambda} $.
\begin{lem}  \label{non3}  Let $ z \in L^\infty_t H^1_\lambda $ and let
$v \in  \tilde{L}^4_{t,\lambda} $ then
\begin{equation}
\| z v \|_{\tilde{L^4}_{t,\lambda}} \lesssim
(\|z\|_{L^\infty_{t,\lambda}} + \|z_x\|_{L^\infty_t L^2_\lambda})
  \|v\|_{\tilde{L}^4_{t,\lambda}} \\
   \quad .
 \label{non3a}
\end{equation}

\end{lem}

\subsection{Proof of \re{nonlinear4}}
 Since $ u $ is real-valued, it holds
$$
\|J_x^s u\|_{L^p_1 L^q_\lambda} \lesssim \|P_1 u\|_{L^p_1 L^q_\lambda}+ \|D_x^s P_{>1} u \|_{L^p_1 L^q_\lambda} \quad .
$$
To estimate the high modes part, we use \re{A4} where we
 expand the exponential function. Hence, we write
\begin{eqnarray}
\| D_x^s P_{>1} u\|_{L^\infty_1 L^2_\lambda}
 & \lesssim & \sum_{k\ge 0} \frac{1}{k!} \|D_x^s( F^k w) \|_{L^\infty_1 L^2_\lambda}
  \nonumber \\
 & &  \hspace*{-10mm} +  \sum_{k\ge 1} \sum_{l\ge 1} \frac{1}{k!\, l!}
  \Bigl\|D_x^s P_{>1}\Bigl( P_{>1}(F^k) \partial_x P_-(F^l)
   \Bigr)\Bigr\|_{L^\infty_1 L^2_\lambda} \; .
   \label{mi60}
\end{eqnarray}
From \re{mi60},  Lemmas \ref{non1} and \ref{non2}, Sobolev
inequalities and  \re{l5}, we infer that for  $ 0\le  s\le 1/2 $,
\begin{eqnarray}
\|D_x^s P_{>1} u\|_{L^\infty_1 L^2_\lambda} &  \lesssim &
 \sum_{k\ge 0} \frac{1}{k!} (\|F^k \|_{L^\infty_{1,\lambda}}
  +  \|\partial_x(F^k) \|_{L^\infty_1 L^2_\lambda} ) \|J_x^s w\|_{L^\infty_1
L^2_\lambda} \nonumber\\
 & &
 \hspace*{-10mm} +  \sum_{k\ge 1} \sum_{l\ge 1} \frac{1}{k!\, l!}
  \|D_x^{5s/4}  P_{>1}(F^k)\|_{L^\infty_1 L^{4/s}_\lambda}
   \|  D_x^{1-s/4} P_-(F^l)\|_{L^\infty_1 L^\frac{4}{2-s}_\lambda} \nonumber\\
    &\lesssim  &
    \sum_{k\ge 0} \frac{1}{k!} (\|F\|_{L^\infty_{1,\lambda}}^k
  + k \|F\|_{L^\infty_{1,\lambda}}^{k-1}\|F_x \|_{L^\infty_1 L^2_\lambda} )
   \| w\|_{Y^s_{1,\lambda}} \nonumber\\
 & &
 \hspace*{-10mm} +  \sum_{k\ge 1} \sum_{l\ge 1} \frac{1}{k!\, l!}
  \|D_x^{s+1/2} \  P_{>1}(F^k)\|_{L^\infty_1 L^2_\lambda}
   \|  \partial_x P_-(F^l)\|_{L^\infty_1 L^2_\lambda} \quad ,
   \label{mi61}
\end{eqnarray}
with
$$
 \Bigl\|D_x^{s+1/2} \  P_{>1}(F^k)\|_{L^\infty_1 L^2_\lambda}
 \lesssim \|\partial_x(F^k)\|_{L^\infty_1 L^2_\lambda}
 \lesssim k
  \|F\|_{L^\infty_{1,\lambda}}^{k-1}\|F_x \|_{L^\infty_1 L^2_\lambda}
  $$
  and
  $$
 \|\partial_x(F^l)\|_{L^\infty_1 L^2_\lambda}
 \lesssim l \,
  \|F\|_{L^\infty_{1,\lambda}}^{l-1}\|F_x \|_{L^\infty_1 L^2_\lambda} \; .
  $$
On the other hand, by the Duhamel formulation of the equation,
the unitarity of $ V(t) $ in $ L^2_\lambda $, the continuity of $ \partial_x P_1 $ in
$ L^2_\lambda $ and Sobolev inequalities, we get
\begin{equation}
\|P_1 u\|_{L^\infty_1 L^2_\lambda}  \lesssim  \|u_0\|_{L^2_\lambda}+
\|u^2\|_{L^1_1 L^2_\lambda}
 \lesssim  \|u_0\|_{L^2_\lambda}+ \|u\|_{L^4_{1,\lambda}}^2\; .
 \label{estP1u}
\end{equation}
This completes the proof of \re{nonlinear4}. \\
\subsection{Proof of \re{nonlinear2}}
\begin{rem}
It would considerably simplify the estimates on $ u $ if we were able to prove that there
 exists $ C >0 $ such that for any $ v\in N_{1,\lambda} $ there exists an extension $ {\tilde
  v} $ of $ v$ satisfying :
  $$
  \|{\tilde v} \|_{X^{7/8,-1}_\lambda}\le C\,  \|{v} \|_{X^{7/8,-1}_{1,\lambda}}
  , \;  \|{\tilde v} \|_{Z^{0,0}_\lambda}\le C\,  \|{v} \|_{Z^{0,0}_{1,\lambda}}
  \mbox{ and }
   \|{\tilde v} \|_{ {\tilde L}^4_{t,\lambda}}\le C\,
\|{v} \|_{{\tilde L}^4_{1,\lambda}}\; .
  $$
  Indeed, we could then take different extensions of $ u$ according to the part of the
$ N_\lambda$-norm we want to estimate. Note, in particular, that
taking the
 extension $ P_{>3} {\tilde u} $ of $ P_{>3} u $ defined by
 $$
 P_{>3} {\tilde u}=\psi(t) \Bigl[ V(t) P_{>3} u_0 +\frac{1}{2}
 \int_0^t V(t-t')P_{>3}\partial_x
 (\psi u )^2(t') \, dt'\Bigr]
 $$
 we directly get
$$
\|P_{>3} u \|_{X^{7/8,-1}_{1,\lambda}}\lesssim
\|u\|_{L^4_{1,\lambda}}^2+ \|u\|_{L^\infty_1 L^2_\lambda}
\lesssim \|u\|_{N_{1,\lambda}}^2+ \|u_0\|_{L^2_\lambda} \; .
$$
\end{rem}
 We start by constructing our
extension $ {F ^*}$ of $ F $.  To construct the high modes part,
we first  need  some how to inverse the map $ F \mapsto W $. From
\re{defW} we infer that
$$
P_{>1}W=P_{>1} (e^{-iF/2})=e^{-iF/2}-P_{\le 1} (e^{-iF/2})
$$
By decomposing $ F $ in $ Q_1 F +P_1 F $, we obtain
$$
e^{-iQ_1 F/2} = e^{iP_1 F/2}\Bigl(P_{>1} W +P_{\le 1}(e^{-iF/2})
\Bigr)
$$
and thus
\begin{eqnarray}
P_{>3} F & = & 2i  P_{>3} \Bigl[ e^{iP_1 F/2}\Bigl(P_{>1} W
-P_{\le
1}(e^{-iF/2})\Bigr)\Bigr]\nonumber \\
 & & -2i   P_{>3}  \Bigl( e^{-iQ_1 F/2} +iQ_1 F/2\Bigr) \quad
 .
 \label{defF}
\end{eqnarray}
Now, let $ {\tilde W} $ be an extension of $ W $ such that $
 \|{\tilde W}_x\|_{M^0_\lambda} \le 2 \| W_x\|_{M^0_{1,\lambda}} $
  and  $ {\tilde F } $ be the extension of $ F $ defined in the
  last section. We set
   \begin{eqnarray}
   P_{>3} F^* & = &  2i \psi \,  P_{>3}\Bigl[ e^{iP_1 \tilde{F}/2}
    \Bigl(  P_{>1} (\psi  {\tilde W})-P_{\le
1}(e^{-i\tilde{F}/2})\Bigr) \Bigr]\nonumber \\
 & & -2i  \psi \, P_{>3}  \Bigl( e^{-iQ_1 \tilde{F}/2} +iQ_1 \tilde{F}/2\Bigr) \quad
 ,
 \label{defFetoile}
\end{eqnarray}
 $ P_{<-3}F^*= \overline{P_{>3}F^*} $ and $ P_3 F^*= P_3 {\tilde
 F}$. It is clear that by construction $ F^*\equiv F $ on $ [0,1] $.
  Note that by \re{lala3}, in Lemma \ref{estimbasfreq}, we already
   have an estimate on the low-modes part $ P_3 F^* $. Moreover,  combining
     estimate \re{lala1} with \re{itit},
    we infer that
    \begin{equation}
\| \tilde{F}\|_{A_\lambda} \lesssim \| \widehat{ \partial_x^{-1}
u_0} \|_{L^1_\xi} + \|u\|^2_{N_{1,\lambda}} + \| u\|_{N_{1,\lambda}}
\quad . \label{K}
    \end{equation}
   To estimate the high-modes part,
     for convenience,
  we drop the $ \; \tilde{}\; $ in the right-hand side of
  \re{defFetoile}. In the remaining of
  this section we assume that $ W $ is supported in time in $
  [-2,2] $.
  \subsubsection{Estimate on the  $ {\tilde L}^4_{t,\lambda}$ -norm}
  Differentiating \re{defFetoile} with respect to $ x$ and
  expanding the exponential function, we get
  \begin{eqnarray*}
  \| P_{>3} F^*_x \|_{{\tilde L}^4_{t,\lambda}} & \lesssim &
  \sum_{k\ge 0} \frac{1}{k!} \Bigl( k \Bigl\|
  (P_1 F_x) (P_1 F)^{k-1} P_{>1} W \Bigr\|_{{\tilde L}^4_{t,\lambda}}
   + \Bigl\| (P_1 F)^{k} P_{>1} W_x \Bigr\|_{{\tilde L}^4_{t,\lambda}}\Bigr) \\
   & & \hspace*{-30mm}+ \sum_{k\ge 0} \sum_{l\ge 0} \frac{1}{k!} \frac{1}{l!}
   \Bigl[  k \Bigl\|\psi \, P_{>3} \Bigl(
  (P_1 F_x) (P_1 F)^{k-1} P_{\le 1} (F^l)\Bigr)  \Bigr\|_{{\tilde L}^4_{t,\lambda}}
  + l \Bigl\| \psi\, P_{>3} \Bigl( (P_1 F)^{k} P_{\le 1} (F^{l-1} F_x)\Bigr)
   \Bigr\|_{{\tilde L}^4_{t,\lambda}}\Bigr]\\
& & + \sum_{k\ge 2} \frac{1}{(k-1)!} \|\psi\,  (Q_1F)^{k-1} Q_1 F_x
 \|_{{\tilde L}^4_{t,\lambda}}
  \end{eqnarray*}
We notice that by the frequency projections,
\begin{equation}
P_{>3} \Bigl( (P_1 F_x) (P_1 F)^{k-1} P_{\le 1} (F^l)\Bigr)
\mbox{ and }
P_{>3} \Bigl( (P_1 F)^k  P_{\le 1} (F^{l-1}F_x)\Bigr) \label{pepe} \quad .
\end{equation}
vanish for $ k\le 2 $. Moreover, decomposing $ P_1 F $ as $ P_1 Q_{\frac{1}{k-1}}F + P_{\frac{1}{k-1}}F $
 we infer that for $ k\ge 3 $ the  two terms appearing in \re{pepe}
are respectively equal to
 $$
 P_{>3}\Bigl[ (P_1 F_x) (P_1 Q_{\frac{1}{k-1}}F)  P_{\le 1} (F^l) G \Bigr]
 \mbox{ and }
  P_{>3}\Bigl[ (P_1 F) (P_1 Q_{\frac{1}{k-1}}F)  P_{\le 1} (F^{l-1} F_x ) G \Bigr]
 $$
  with
  $$
G= \sum_{q=1}^{k-1} C_{k-1}^q   (P_1 Q_{\frac{1}{k-1}}F)^{q-1}
(P_{\frac{1}{k-1}}F)^{k-1-q}\quad .
$$
Note that $ G $ can be also written as
$$
G= \sum_{j=0}^{k-2} \frac{C_{k-1}^{j+1}}{C_{k-2}^j} C_{k-2}^j
\, (P_1 Q_{\frac{1}{k-1}}F)^{j} (P_{\frac{1}{k-1}}F)^{k-2-j}
$$
and thus it is not too hard to see that
\begin{equation}
\|G\|_{L^\infty_{t,\lambda}}\lesssim
\|G \|_{A_\lambda} \lesssim (k-1) \|F\|_{A_\lambda}^{k-2} \quad .\label{pepe2}
\end{equation}
Therefore,
 using that, by Sobolev inequalities,
$$
\| Q_{\frac{1}{k-1}} P_1 F\|_{L^\infty_{t,\lambda}} \lesssim
(k-1) \|F_x \|_{L^\infty_t L^2_\lambda}\quad ,
$$
  using  Lemma \ref{non3}  and the embedding
 $ X_\lambda^{1/2,0} \hookrightarrow {\tilde L}^4_{2,\lambda} $ (see \re{l4}), we infer that
\begin{eqnarray}
  \| P_{>3} F^* \|_{{\tilde L}^4_{t,\lambda}} & \lesssim &\Bigl[(
 \|F_x\|_{L^\infty_t L^2_\lambda} +1)
 \| W_x\|_{{\tilde L}^4_{t,\lambda}} \nonumber \\
  & &+ \|F_x\|_{L^\infty_t L^2_\lambda}
  \|\psi_2  F_x\|_{{\tilde L}^4_{t,\lambda}}
  \Bigr]
  e^{2\|F\|_{L^\infty_{t,\lambda}}}\nonumber\\
  & \lesssim & \Bigl( \|w\|_{M^0_{1,\lambda}} ( 1+ \|u\|_{N_{1,\lambda}})
  + \|u\|_{N_{1,\lambda}}^2 \Bigr)\, e^{\tilde K}  \quad .\label{co2}
  \end{eqnarray}
  where $ {\tilde K} $ is defined as in \re{nonlinear3}.
  \subsubsection{Estimate on the $ Z^{0,0}_\lambda$ -norm}
Now, using again the frequency projections and that $ A_\lambda $
is clearly an algebra, we deduce from \re{defFetoile} and \re{pepe2} that
\begin{eqnarray*}
\| P_{>3} F^*_x \|_{Z^{0,0}_\lambda} & \lesssim & \sum_{k\ge 0}
\frac{1}{k!} \|P_1 F\|_{A_\lambda}^{k-1}\Bigl( k \|P_1 F_x\|_{Z^{0,0}_\lambda}
\|P_{>1} W\|_{A_\lambda}+\|P_1 F\|_{A_\lambda} \|W_x\|_{Z^{0,0}_\lambda} \Bigr) \nonumber  \\
 & & +  \sum_{k\ge 3}  \sum_{l\ge 0}\frac{k}{k!\, l!}(k+l)
 \|Q_{\frac{1}{k-1}} P_1 F\|_{A_\lambda} \|F_x\|_{Z^{0,0}_\lambda} \| F\|_{A_\lambda}^{k+l-2}\\
 &  &
 +\sum_{k\ge 2} \frac{k}{k!} \|F_x\|_{Z^{0,0}_\lambda
 }\|Q_1  F\|_{A_\lambda}^{k-1}\quad .
\end{eqnarray*}
Using that, by Cauchy-Schwarz in $ \xi $,
$$
\|P_{>1} W\|_{A_\lambda} \lesssim \|W_x\|_{Z^{0,0}_\lambda}\mbox{ and }\;
\|\widehat{Q_{\frac{1}{k-1}} P_1 F}\|_{A_\lambda}\lesssim \,
(k-1) \, \|  F_x\|_{Z^{0,0}_\lambda}
$$
we infer that
\begin{eqnarray}
  \| P_{>3} F^*_x \|_{Z^{0,0}_{\lambda}} & \lesssim &\Bigl(
 \|F_x\|_{N_\lambda}  \| W_x\|_{Z^{0,0}_\lambda}
 + \| W_x\|_{Z^{0,0}_{\lambda}}+ \|F_x\|_{N_\lambda}^2\Bigr)
  e^{2\|F\|_{A_\lambda}} \nonumber\\
  & \lesssim & \Bigl( \|w\|_{M^0_{1,\lambda}} ( 1+ \|u\|_{N_{1,\lambda}})
  + \|u\|_{N_{1,\lambda}}^2 \Bigr)\, e^{\tilde K}  \quad .\label{co3}
  \end{eqnarray}

\subsubsection{ Estimate on the $ X^{7/8,-1}_\lambda $-norm}
It remains to estimate the $ X^{7/8,-1}_\lambda $-norm of $ P_{>3}
F^*_x $. Note that obviously
$$
\| P_{>3} F^*_x \|_{X^{7/8,-1}_\lambda} \sim \| P_{>3} F^* \|_{X^{7/8,0}_\lambda}
$$
From \re{defFetoile} we infer  that
\begin{eqnarray*}
\| P_{>3} F^* \|_{X^{7/8,0}_\lambda} & \lesssim & \sum_{k\ge 0}
\frac{1}{k!} \Bigl\|P_{>3}\Bigl(  (P_1 F)^k P_{>1} W\Bigr) \Bigr\|_{X^{7/8,0}_\lambda} \\
 & & + \sum_{k\ge 3}  \sum_{l\ge 0}\frac{1}{k!\, l!}
 \Bigl\|\psi \, P_{>3} \Bigl(  (P_1 F)^{k}
  P_{\le 1} (F^l)\Bigr) \Bigr\|_{X^{7/8,0}_\lambda} \\
 & & +\sum_{k\ge 2} \frac{1}{k!}\Bigl\|\psi \, P_3 \Bigl( Q_1  F)^k\Bigr)
 \Bigr\|_{X^{7/8,0}_\lambda} \\
  & =& \sum_{k\ge 0}\,  \frac{1}{k!} \, I_k + \sum_{k\ge 3} \sum_{l\ge 0}\, \frac{1}{k!\,
  l!}\, J_{k,l}+\sum_{k\ge 2} \, \frac{1}{k!} \, L_k \quad .
\end{eqnarray*}
Let us  estimate $ I_k $, $ J_{k,l}$  and $ L_k $, one by one.\vspace*{2mm} \\
{\bf  i) Estimate on $ I_k $.} First note that for $ k=0 $, we
have directly
\begin{equation}
I_0 \lesssim \|P_{>1} W \|_{X^{7/8,0}_\lambda} \lesssim \|w
\|_{X^{7/8,-1}_\lambda}\lesssim \|w\|_{M^0_{1,\lambda}} \quad .\label{yuyu1}
\end{equation}
Now, for $ k\ge 1 $,
$$I_k=\Bigl\| \chi_{\{\xi\ge 3\}} \langle
\sigma\rangle^{7/8}\int_{\R^{k}\times (\lambda^{-1}\Z)^{k}}
\widehat{P_1 F}(\tau_1,\xi_1) ..
 \widehat{P_1 F}(\tau_k, \xi_k) \, \widehat{P_{>1} W}(\tau_{k+1},\xi_{k+1}) \Bigr\|_{L^2_{\tau,\xi}}
$$
where $ \sigma=\tau+\xi|\xi| $ and $  \Bigl(\sum_{i=1}^{k+1}
\tau_i, \sum_{i=1}^{k+1} \xi_i\Bigr)=(\tau,\xi) $ . \\
We divide $\R^{k+1} \times (\lambda^{-1}\Z)^{k+1} $ in different regions. \vspace{2mm} \\
$ \bullet $ The region $ |\sigma|\le 2^{10} k  $. In this region,
clearly,
\begin{equation}
I_k \lesssim k \|P_{>1} W \|_{L^\infty_t
L^2_\lambda}\|F\|_{A_\lambda}^k  \lesssim k
\|w\|_{M^0_{1,\lambda}} \|F\|_{A_\lambda}^k\quad . \label{yuyu2}
\end{equation}
$ \bullet $ The region \{$  2^{4} k \,
|\tau_{k+1}+\xi_{k+1}|\xi_{k+1}||\ge |\sigma|\mbox{ and }
|\sigma|>2^{10}k \} $.
 In this region it is easy to see that
\begin{eqnarray}
I_k \lesssim k \|P_{>1} W \|_{X^{7/8,0}_\lambda}
\|F\|_{A_\lambda}^k \lesssim k \|Q_1 w\|_{X^{7/8,-1}_\lambda}
\|F\|_{A_\lambda}^k\lesssim k \|w\|_{M^0_{1,\lambda}}
\|F\|_{A_\lambda}^k \quad . \label{yuyu3}
\end{eqnarray}
$\bullet $  The region $\{\exists i \in  \{1,..,k\} , \quad 2^4 k \,
|\tau_i+\xi_i|\xi_i||\ge  \langle\sigma\rangle\,
 \mbox{ and }
|\sigma|>2^{10}k\}$. Then we have
\begin{eqnarray}
I_k & \lesssim & k \, \|P_1 F\|_{{\dot X}^{7/8,0}_\lambda}
 \|P_1 F\|^{k-1}_{A_\lambda}
 \|P_{>1} W\|_{A_\lambda} \nonumber \\
& \lesssim & k \, \|P_1 F\|_{{\dot X}^{1,0}_\lambda}
 \|P_1 F\|^{k-1}_{A_\lambda}
 \|P_{>1} W_x \|_{Z^{0,0}_\lambda} \nonumber \\
 & \lesssim & k \|w\|_{M^0_{1,\lambda}}\|P_1 F\|_{{\dot X}^{1,0}_\lambda}
\|P_1 F\|^{k-1}_{A_\lambda} \nonumber\\
& \lesssim & k\,
\|w\|_{M^0_{1,\lambda}}(\|u\|_{N_{1,\lambda}}+\|u\|_{N_{1,\lambda}}^2)
\| F\|^{k-1}_{A_\lambda} \label{yuyu4}
\end{eqnarray}
where we used \re{lala2} in the last step.\\
 $ \bullet $ The region
 $ \{ |\sigma| \ge 2^4 k  \displaystyle\max_{i=1,..,k+1}
 |\tau_i+\xi_i|\xi_i|| \mbox{ and }  |\sigma|>2^{10} k \}$.
In this region, since $ \xi\ge  0 $,   we have
\begin{equation}
\langle \sigma \rangle \le 2    |\sigma| \le 4  \Bigl|\sigma
-\sum_{i=1}^{k+1} (\tau_i-\xi_i|\xi_i|)\Bigr| =
\Bigl|\Bigl(\sum_{i=1}^{k+1} \xi_i \Bigr)^2-\sum_{i=1}^{k+1}
\xi_i |\xi_{i}|  \Bigr| \ \label{c17}\; .
\end{equation}
Let us denote by $ |\xi_{i_1}  |= \max|\xi_i| $ and $
\displaystyle|\xi_{i_2}|=\max_{i\neq i_1} |\xi_i| $. We claim
that \re{c17} implies
\begin{equation}
\langle \sigma \rangle \le  2^5 k^2\, |\xi_{i_1}  |\,
|\xi_{i_2} |\quad . \label{sigsig}
\end{equation}
Indeed, either $ 2 k |\xi_{i_2}| \ge |\xi_{i_1}| $ and
$$
 \langle\sigma\rangle \le 4 \Bigl(\Bigl(\sum_{i=1}^{k+1} |\xi_i| \Bigr)^2+\sum_{i=1}^{k+1}
 |\xi_{i}|^2  \Bigr) \le 2^5 k^2\,
|\xi_{i_1}  |\, |\xi_{i_2} |\quad
$$
or $  |\xi_{i_1}| \ge 2  k |\xi_{i_2}| $ and then $ \xi $ and
$\xi_{i_1} $ have the same sign so that
$$
\langle \sigma \rangle \le 4 \Bigl( \sum_{i\neq i_1} |\xi_1|^2 +
\Bigl( \sum_{i\neq i_1} |\xi_i|\Bigr)^2 +2 |\xi_{i_1}|
 \sum_{i\neq i_1} |\xi_i |\Bigr)  \le 2^4
k^2  |\xi_{i_1}  |\, |\xi_{i_2} |\quad .
$$
From \re{sigsig}, we infer that in this region,
\begin{eqnarray}
I_k & \lesssim & k \|P_1 F_x \|_{A_\lambda} \|P_{>1}
W_x\|_{L^2_{t,\lambda}} \|F\|^{k-1}_{A_\lambda} + k(k-1) \|P_1 F_x
\|_{A_\lambda}^2 \|P_{>1} W\|_{L^2_{t,\lambda}}
\|F\|^{k-2}_{A_\lambda} \nonumber \\
 &\lesssim & k^2 \|P_{>1} W_x \|_{L^2_{t,\lambda}}
\| F\|_{A_\lambda}^{k-1} \nonumber\\
  &\lesssim & k^2 \|w\|_{M^0_{1,\lambda}}\, \| F\|_{A_\lambda}^{k-1}\quad .
  \label{yuyu5}
\end{eqnarray}
 {\bf  ii) Estimate on $ J_{k,l}$.}
 Proceeding as in the treatment of the terms in
 \re{pepe}, we can write $ J_{k,l} $ as
$$
J_{k,l}=\Bigl\|P_{>3} \Bigl( (\psi_2 P_1 Q_{\frac{2}{k-1}} F)^2  P_{\le 1} (F^l)\, G   \Bigr\|_{X^{7/8,0}_\lambda}
$$
where
$$
G= \sum_{q=2}^k C^q_k (P_1 Q_{\frac{2}{k-1}} F)^{q-2}
(P_{\frac{2}{k-1}}F)^{k-q}=\sum_{j=0}^{k-2}
\frac{C^{j+2}_k}{C^j_{k-2}} C_{k-2}^j (P_1 Q_{\frac{2}{k-1}}
F)^{j} (P_{\frac{2}{k-1}}F)^{k-2-j}\quad .
$$
Clearly
\begin{equation}
|\hat{G}|\lesssim k(k-1) |\widehat{P_1 F}|\ast \cdot \cdot |\widehat{P_1 F}|\label{pepe3}
\end{equation}
and thus
\begin{equation}
 \|G\|_{A_\lambda} \lesssim k^2 \|P_1 F\|_{A_\lambda}^{k-2} \; .\label{pepe4}
\end{equation}
 Note first that we can assume that $ |\sigma|\ge 2^{10}(k+l) $ since otherwise obviously,
 \begin{eqnarray}
J_{k,l} & \lesssim &  k^2\|\psi_2 Q_\frac{2}{k-1} P_1
F\|_{L^4_{t,\lambda}}^2 \|F\|_{A_\lambda}^{k+l-2} \lesssim k^4
\|\psi_2 F_x\|_{L^4_{t,\lambda}}^2 \|F\|_{A_\lambda}^{k+l-2}
\nonumber \\
&  \lesssim & k^4 \|u\|_{N_{1,\lambda}}^2
\|F\|_{A_\lambda}^{k+l-2}\quad . \label{yuyu6}
\end{eqnarray}
We have thus to estimate
 \begin{eqnarray*}
 {\tilde J}_{k,l} & = & \Bigl\|
\chi_{\{\xi\ge 3\}}\chi_{\{|\sigma|\ge 2^{10}k\}} \langle
\sigma\rangle^{7/8}{\cal F}_{t,x} \Bigl((\psi_2 Q_\frac{2}{k-1} P_1 F)^2
  P_{\le 1} (F^l)\, G \Bigr)(\tau,\xi)  \Bigr\|_{L^2_{\tau,\xi}}
\end{eqnarray*}
where $ \sigma=\tau+\xi|\xi| $. \\
As in Lemma \ref{nonlin2}, one of  the difficulties is that we do not know if
 $ {\cal F}^{-1}_{t,x}(|\hat{F}|) $ belongs to $ L^4_{t,\lambda} $. Using again
  the Littlewood-Paley decomposition it can be seen that for $ l\ge 2 $,
\begin{equation}
F^l = \sum_{i_1\ge i_2\ge 0} \hspace*{-2mm}\Delta_{i_1}(F) \Delta_{i_2}(F) \sum_{0\le
 i_3,..,i_{l}\le  i_2} \hspace*{-5mm} n(i_1,..,i_l)\prod_{j=3}^{l} \Delta_{i_j} (F)
\; ,\label{dec}
\end{equation}
where  $ n(i_1,..,i_l) $ is an integer belonging to
 $\{1,..,l(l-1)\} $ (Note for instance that $ n(i_1,..,i_l)=1 $ for $ i_1=\cdot\cdot= i_{l} $
   and $ n(i_1,..,i_l)=l(l-1) $ for $ i_1\neq\cdot\cdot\neq i_{l} $).
We set
$$
H_{j,q,l}= \Delta_j(F) \Delta_q(F) \sum_{0\le
 i_3,..,i_{l}\le  q} \hspace*{-5mm} n(j,q,i_3,..,i_l)\prod_{m=3}^{l}
\Delta_{i_m} (F) \quad .
 $$
 It is clear that for $ l\ge 2 $,
 \begin{eqnarray*}
 {\tilde J}_{k,l} & \lesssim & \sum_{j\ge q\ge 1} \Bigl\|
 \chi_{\tiny \{\begin{array}{l} \xi\ge 3 \\|\sigma|\ge 2^{10}k \end{array}\!\! \}}
\langle
\sigma\rangle^{7/8}{\cal F}_{t,x} \Bigl((\psi_2 Q_\frac{2}{k-1} P_1 F)^2 \, G \,
  P_{\le 1} (H_{j,q,l})\Bigr)
  \Bigr\|_{L^2_{\tau,\xi}} \\
  & & +l
 \Bigl\|
 \chi_{\tiny \{\begin{array}{l} \xi\ge 3 \\|\sigma|\ge 2^{10}k \end{array}\!\! \}}
\langle \sigma\rangle^{7/8} \sum_{j\ge 0}{\cal F}_{t,x}
\Bigl((\psi_2 Q_\frac{2}{k-1} P_1 F)^2 \, G \,
  P_{\le 1}\Bigl(\Delta_j(F) \Delta_0(F)^{l-1} \Bigr)\Bigr\|_{L^2_{\tau,\xi}}\\
  & =&  \Lambda_{k,l}+\Gamma_{k,l} \quad .
\end{eqnarray*}
Let us write $ \Lambda_{k,l} $ as the sum of two terms :
\begin{eqnarray*}
 \Lambda_{k,l} & = &
 \sum_{j\ge q\ge 1} \Bigl\|
\chi_{\{\xi\ge 3\}}\chi_{\{ |\sigma|\in D^1_k\}} \langle
\sigma\rangle^{7/8}{\cal F}_{t,x} \Bigl((\psi_2 Q_\frac{2}{k-1}
P_1 F)^2 G
  P_{\le 1} (H_{j,q,l})\Bigr) \Bigr\|_{L^2_{\tau,\xi}}\\
 &  & +  \sum_{j\ge q\ge 1} \Bigl\|
  \chi_{ \{\xi\ge 3 , \,
  |\sigma|\in D^2_k\}}
 \langle
\sigma\rangle^{7/8}{\cal F}_{t,x} \Bigl((\psi_2 Q_\frac{2}{k-1} P_1 F)^2 G
  P_{\le 1} ((H_{j,q,l})\Bigr)  \Bigr\|_{L^2_{\tau,\xi}} \\
  & = & \Lambda_{k,l}^1+\Lambda_{k,l}^2 \quad ,
\end{eqnarray*}
with  $$ D^1_k=[2^{10}(k+l),(k+l)^2 2^{8+j+q}[ \mbox{ and }
  D^2_k=[\max\Bigl((2^{10}(k+l) ,(k+l)^2 2^{8+j+q}\Bigr), +\infty[ \quad .
$$
 From the definition of $ H_{j,q,l} $, \re{pepe4}  and \re{estA}
we infer that for $ l\ge 2 $,
\begin{eqnarray}
\Lambda_{k,l}^1 & \lesssim & k^2 \sum_{j=1}^\infty \sum_{q=1}^{j}
2^{7j/8} 2^{7q/8} \,
 \|F\|_{A_\lambda}^k
  \|\psi_2^2 \, H_{j,q,l}\|_{L^2_{\tau,\xi}} \nonumber \\
  & \lesssim & (kl)^2
 \|F\|_{A_\lambda}^{k+l-2}
   \Bigl(\sum_{j=1}^\infty  2^{7j/8} \,\|\psi_2 \Delta_j F \|_{L^4_{t,\lambda}}\Bigr)^2
    \nonumber \\
  & \lesssim & (kl)^2   \|F\|_{A_\lambda}^{k+l-2} \|\psi_2
  F_x\|_{L^4_{t,\lambda}}^2
   \lesssim  (kl)^2\|F\|_{A_\lambda}^{k+l-2} \|F_x\|_{N_\lambda}^2\quad . \label{didi}
\end{eqnarray}
On the other hand , using \re{pepe3}, it is easy to check  that for $ l\ge 2 $,
\begin{eqnarray*}
\Lambda_{k,l}^2 & \lesssim & (kl)^2 \Bigl\| \chi_{\{\xi\ge 3\}}
\langle \sigma\rangle^{7/8}\int_{\R^{k+l-1}\times
(\lambda^{-1}\Z)^{k+l-1}} \chi_{\{|\sigma| > 2^5(k+l)^2
|\xi_{i_1}\xi_{i_2}|\}}|\widehat{Q_\frac{2}{k-1} P_1
F}(\tau_1,\xi_1)| \\
& & |{\cal F}_{t,x}\Bigl( \psi_2 Q_\frac{2}{k-1} P_1 F\Bigr)(\tau_2,\xi_2) |\, |\widehat{
P_1 F}(\tau_3,\xi_3)| |\widehat {P_1 F}(\tau_k,\xi_k)|
 |\widehat{ F}(\tau_{k+1}, \xi_{k+1})| .. |\widehat{F}(\tau_{k+l},\xi_{k+l})| \Bigr\|_{L^2_{\tau,\xi}}
 \end{eqnarray*}
 where $ |\xi_{i_1}  |= \max|\xi_i| \mbox{ and }
|\xi_{i_2}|=\max_{i\neq i_1} |\xi_i|\;$. But the same
considerations as in
 \re{c17}-\re{sigsig} ensure that $2^{10} (k+l)\le  |\sigma |\le 10 (k+l)
 \max_{i=1,..,k+l} |
  \tau_i-\xi_i|\xi_i|| $ in the region of integration above. Therefore,
   according to Lemma \ref{estimbasfreq},
\begin{eqnarray}
\Lambda_{k,l}^2 & \lesssim & (kl)^2 (k+l)   \| \psi_2 Q_\frac{2}{k-1}
P_1 F\|_{A_\lambda}^2
\|F\|_{A_\lambda}^{k+l-3} \|F\|_{{\dot X}^{7/8,0}_{\lambda}} \nonumber \\
& & +  (kl)^2(k+l)\| \psi_2 Q_\frac{2}{k-1} P_1 F\|_{A_\lambda}
\|F\|_{A_\lambda}^{k+l-2} \|\psi_2 Q_\frac{2}{k-1} P_1 F\|_{X^{7/8,0}_{\lambda}}
\nonumber  \\
 &\lesssim &  (kl)^2 (k+l)k  \|F_x\|_{Z^{0,0}_\lambda} \|F\|_{A_\lambda}^{k+l-3}\nonumber\\
 & & \Bigl( k \|F_x\|_{Z^{0,0}_\lambda} ( \|P_3 F\|_{{\dot X}^{1,0}_\lambda} + \|Q_3 F\|_{X^{7/8,0}_\lambda})
+ \|F\|_{A_\lambda}
   \|\psi_2 Q_\frac{2}{k-1} P_1 F\|_{X^{7/8,0}_{\lambda}} \Bigr)\nonumber \\
  & \lesssim & (k+l)^5 \|u\|_{N_{1,\lambda}}^2
  (\|u\|_{N_{1,\lambda}}+\|F\|_{A_\lambda})(1
  + \|u\|_{N_{1,\lambda}} )\|F\|_{A_\lambda}^{k+l-3} \label{fi1}\, .
 \end{eqnarray}
 It remains to estimate $ \Gamma_{k,l} $ for $ l\ge 2 $.
We notice that
\begin{eqnarray*}
\Gamma_{k,l} & \lesssim & k^2 l \Bigl\|\int_{\R^{k+l-1}\times
(\lambda^{-1} \Z)^{k+l-1}} |{\cal F}_{t,x}\Bigl( \psi_2
Q_\frac{2}{k-1} P_1 F\Bigr)(\tau_1,\xi_1)|
|{\cal F}_{t,x}\Bigl(\psi_2 Q_\frac{2}{k-1} P_1 F\Bigr)(\tau_2,\xi_2)| \\
& & \quad \quad \displaystyle\prod_{i=3}^{k} |\widehat{P_1
F}(\tau_i,\xi_i)|\; |\hat{F}(\tau_{k+1}, \xi_{l+1})|\;
\prod_{i=k+2}^{k+l} |\widehat{P_3 F}(\tau_i,\xi_i)|
\Bigr\|_{L^2_{\tau,\xi}}
\end{eqnarray*}
which can be estimated in the same way we did for $ I_k $.
 More precisely,
in the region,
  $ 2^4 (k+l) \displaystyle\max_{i=1,..,k+l} |\tau_i-\xi_i|\xi_i||
   \ge |\sigma|  $ we easily get as above
\begin{equation}
{\Gamma}_{k,l}  \lesssim (k+l)^5 \|u\|_{N_{1,\lambda}}^2
(\|u\|_{N_{1,\lambda}}+\|F\|_{A_\lambda})(1
  + \|u\|_{N_{1,\lambda}} )
   \Bigr)\|F\|_{A_\lambda}^{k+l-3},   \label{fi2}
\end{equation}
and in the region $\displaystyle |\sigma | \ge 2^4 (k+l)
\max_{i=1,..,k+l}| \tau_i-\xi_i|\xi_i||  $ we infer from
\re{sigsig} that
\begin{equation}
{\tilde J}_{k,l}  \lesssim (k+l)^5 \|u\|_{N_{1,\lambda}}^2
\|F\|_{A_\lambda}^{k+l-2}\;  .
 \label{fi3}
\end{equation}
 Finally, we notice that  $ {\tilde J}_{k,0}$ and ${\tilde J}_{k,1} $
 with $ k\ge 3 $ can be estimated exactly in the same way.\\
{\bf  iii) Estimate on  $L_k $ } This term can be treated in the same way
as the preceding one
 and   is even much simpler. Since $ k\ge 2 $ we can decompose $ Q_1(F)^k $ as
  we did
  for $F^l $ in \re{dec} and then proceed exactly in the same way as for
  $ J_{k,l} $. We get
\begin{equation}
 L_{k}\lesssim   k^5  \|u\|_{N_{1,\lambda}}^{k}(1+  \|u\|_{N_{1,\lambda}}^{k})\vspace*{2mm} \label{yuyu7}
 \vspace*{2mm}
 \end{equation}
 Gathering \re{yuyu1}-\re{yuyu4}, \re{yuyu5}-\re{yuyu6} and \re{didi}-\re{yuyu7}, we finally deduce that
 \begin{equation}
\|P_{>3} F^*_x\|_{X^{7/8,-1}_\lambda} \lesssim \Bigl(
\|w\|_{M^0_{1,\lambda}}
(1+\|u\|_{N_{1,\lambda}})+\|u\|_{N_{1,\lambda}}^2 \Bigr) \, e^{\tilde K}
\quad \label{Xnorm}
 \end{equation}
which ends the proof of \re{nonlinear2}.
\section{ Uniform estimates and Lipschitz bound for small initial data}
\label{section6}
  \subsection{Uniform estimate for  small initial data}
  We are now ready to state the following crucial proposition on the uniform boundedness
   of small smooth solutions to (BO).
 \begin{pro} \label{propuniformbound}
 Let $ 0\le s\le 1/2 $ and   $ K\ge 1 $ be given.  There exists
 $ 0<\varepsilon:=\varepsilon(K)\sim e^{-8CK} <1 $ such that for any $ u_0\in {H}^\infty_{0,\lambda} $
  with
  $$
\|\widehat{\partial_x^{-1} u_0}\|_{L^1_\xi} \lesssim K \quad
\mbox{ and } \quad \|u_0\|_{L^2_\lambda} \lesssim \varepsilon^2
\quad ,
  $$
  the emanating solution $u\in C(\R;{H}^\infty_{0,\lambda})$ to (BO)
 satisfies
 \begin{equation}
 \|u\|_{L^\infty_1 H^s_\lambda} \lesssim e^{2C K} \|u_0\|_{H^s_\lambda}\quad\mbox{ and }
 \quad \|w\|_{M^s_{1,\lambda}}\lesssim
  e^{K}\|u_0\|_{H^s_\lambda}\quad .
 \label{ets}
 \end{equation}
 \end{pro}
 {\it Proof.}
  For $ K\ge 1 $ given, let  $ B_{K,\lambda} $ be  the  small closed ball  of $ {L}^2_{\lambda} $
  defined by
\begin{equation}
B_{K,\lambda}  := \Bigl\{ \varphi\in { L}^2_{\lambda}, \, \;
\|\widehat{\partial_x^{-1} \varphi} \|_{L^1_\xi} \lesssim K
\mbox{ and } \|\varphi \|_{L^2_\lambda} \lesssim \varepsilon(K)^2
\, \Bigr\} \label{o1}
\end{equation}
 where $ 0<\varepsilon(K)\sim e^{-8C K} <\! \! < 1 $ ($ C>1$ is the universal constant
  appearing in \re{nonlinear3})  only depends on $ K $
  and the implicit constants contained in the estimates of the preceding sections. At this stage, it worth recalling that
 these  implicit constants do not depend on the period $ \lambda $. \\
 We set  $\varepsilon:=\varepsilon(K) $. For
  $ u_0$ belonging to  $ {H}^\infty_{0,\lambda} \cap B_{K,\lambda}
  $,
   we want    to  show that the emanating solution
   $ u\in C(\R;{H}^\infty_{0,\lambda}) $, given by the classical well-posedness results
    (cf. \cite{ABFS}, \cite{Io}),
    satisfies
\begin{equation}
    \| u\|_{N_{1,\lambda}}\lesssim e^{2 C K} \varepsilon^2 \quad \mbox{ and }
    \quad  \| w\|_{M^{0}_{1,\lambda}}\lesssim e^{ K} \varepsilon^2 \quad . \label{o2}
    \end{equation}
     \re{ets} then obviously follows
    from
   \re{o2} together with \re{nonlinear1} and \re{nonlinear4}.\\
     Clearly, since $ u $ satisfies the equation, $ u $ belongs in fact to $ C^\infty(\R;H^\infty_\lambda) $ and thus $   u $ and $ w  $ belong to $ {M^{\infty}_{1,\lambda}}\cap N_{1,\lambda} $.  We are going  to implement a  bootstrap  argument.
      Since we have chosen to take $ T=1 $ we can not use any continuity argument in time
  but  as in \cite{PlaBu} we will apply a continuity argument on the space period. Recall  that if $ u(t,x) $ is a $2\lambda \pi $-periodic solution of (BO) on $[0,T] $ with initial data
 $ u_0 $ then $ u_\beta(t,x)=\beta^{-1} u(\beta^{-2}t,\beta^{-1} x)$ is a $ (2\pi \lambda\beta) $-periodic solution of (BO) on
 $ [0,\beta^2 T] $
 emanating from $u_{0,\beta}=\beta^{-1} u_0(\beta^{-1} x) $. Moreover,
   denoting by $ w_\beta $ the gauge transform of $ u_\beta $, it is worth noticing that
 \begin{equation}
  w_{\beta}(t,x)=\beta^{-1} w(\beta^{-2}t,\beta^{-1} x) \label{wbeta} \quad .
 \end{equation}
 Straightforward  computations give
 \begin{equation}
 \|u_{0,\beta} \|_{L^2_{\lambda \beta}} = \beta^{-1/2} \|u_0\|_{L^2_{\lambda}}
 \mbox{ and } \|\widehat{\partial_x^{-1} u_{0,\beta}} \|_{L^1_{\xi}} =
  \|\widehat{\partial_x^{-1}u_0}\|_{L^1_{\xi}} \; . \label{init1}
 \end{equation}
Note that $\|\widehat{\partial_x^{-1} u_{0}}
 \|_{L^1_{\xi}}$ is  invariant by the symmetry dilation of
 (BO). In  the same way  one can easily check that the $
N_{1,\lambda\beta} $-norm  of
   $ u_\beta $ and the $ M^0_{1,\lambda\beta}$ norm of $w_\beta $    tend
    to  $ 0 $ as $ \beta $ tends to infinity.
  Hence, for $ \beta $ large enough,  $u_\beta $ and $ w_\beta $       satisfy
 \begin{equation}
 \| u_\beta \|_{N_{1,\lambda\beta}}
     + \|  w_\beta \|_{M^{0}_{1,\lambda\beta }} \lesssim
  \varepsilon\; . \label{init2}
\end{equation}
  \re{nonlinear1}  then clearly ensures
that $ \| w_\beta\|_{M^{0}_{1,\lambda\beta }} \lesssim
(1+\|u_{0,\beta}\|_{L^2_{\lambda\beta }}) e^{K}
\|u_{0,\beta}\|_{L^2_{\lambda\beta}}$
 and \re{nonlinear2}-\re{nonlinear3}   ensure that
 $$
  \|u_\beta \|_{N_{1,\lambda\beta }} \lesssim  (1+\|u_{0,\beta}\|_{L^2_{\lambda\beta}
  })e^{ 2 C K}
    \|u_{0,\beta}\|_{L^2_{\lambda\beta} } \quad .
  $$
 Therefore, by the assumptions on $ u_0 $ and \re{init1}, we finally get
  \begin{equation}
   \| u_\beta \|_{N_{1,\lambda\beta}}\lesssim  e^{2 C K} \beta^{-1/2} \varepsilon^2
     \mbox{ and }
     \|  w_\beta \|_{M^{0}_{1,\lambda\beta }} \lesssim e^{K} \beta^{-1/2} \varepsilon^2
    \quad
     \label{init3}
  \end{equation}
  which,  by the definition of  $\varepsilon $, proves that
  $$
\| u_\beta \|_{N_{1,\lambda\beta}}+ \|  w_\beta \|_{M^{0}_{1,\lambda\beta }} \lesssim
\beta^{-1/2} \varepsilon^{3/2} \quad .
  $$
  $ \beta\mapsto  \| u_\beta \|_{N_{1,\lambda\beta}}
     +\|  w_\beta \|_{M^{0}_{1,\lambda\beta }} $ being clearly
     continuous, a
     classical
    continuity argument in $\beta $ ensures that  we can take $ \beta=1
    $ in \re{init3}. This completes the proof of    \re{o2} and thus of
    \re{ets}.
\subsection{Lipschitz bound}
To prove the continuity of the solution as well as the continuity
 the flow-map we will derive  a Lipschitz bound  on the solution-map $ u_0 \mapsto u $
 for small solutions of $(BO) $ (Note that up to now this map in only defined
 on $ H^\infty_\lambda $).\vspace*{2mm} \\
Let  $ u_1 $ and $ u_2 $ be two solutions of (BO) in $
N_{1,\lambda}\cap C([0,T];H^s_\lambda) $ associated with initial
data $ \varphi_1 $ and $ \varphi_2 $
 in $ B_{K,\lambda}\cap {H}^{s}_{\lambda} $ such that their gauge transforms
 $ w_1 $ and $ w_2$ belong to $ M^s_{1,\lambda} $.
 We assume that they satisfy
\begin{equation}
   \|u_i\|_{N_{1,\lambda}}
   + \| w_i\|_{M^{0}_{1,\lambda}}\lesssim \varepsilon^2 , \quad i=1,2\quad , \label{oo2}
    \end{equation}
   where $ 0<\varepsilon=\varepsilon(K)<\!\! <1$.\\
We set  $ W_i=P_+(e^{-iF_i/2}) $ with $ F_i=\partial_x^{-1} u_i $,
 $ w_i=\partial_x W_i$,  $ v=u_1-u_2 $, $Z=W_1-W_2 $ and $ z=Z_x $. \\
It is easy to check that
\begin{eqnarray}
v
 & = & 2 i e^{iF_1/2} \Bigl[ z +\partial_x P_- \Bigl( e^{-iF_1/2}  -e^{-iF_2/2}
 \Bigr)\Bigr] \nonumber \\
 & &  +2i ( e^{iF_1/2} -e^{iF_2/2}) \Bigl(w_2
  +  \partial_x P_- (e^{-iF_2/2}) \Bigr) \label{o7}
\end{eqnarray}
and that $ z $ satifies
 \begin{eqnarray}
 z_t-i z_{xx}  & = & -\partial_x P_+ \Bigl[ W_1
 \partial_x P_-(v) \Bigr]
 -\partial_x P_+ \Bigl[ Z \, P_-(\partial_x u_2)\Bigr] \nonumber \\
  & & +\frac{i}{4} \Bigl( P_0(u_1^2) z+ P_0(u_1^2-u_2^2) w_2\Bigr)\; . \label{o13}
 \end{eqnarray}
 As in the obtention of \re{eq2w}, we substitute \re{o7} in \re{o13} to get
 \arraycolsep1pt
 \begin{eqnarray*}
 z_t-i z_{xx}  & = & 2i \partial_x P_+ \Bigl[W_1
 \partial_x P_-(e^{-iF_1/2} \overline{z} +
  (e^{-iF_1/2}-e^{-iF_1/2}) \overline{w_2}) \Bigr] \nonumber \\
  & & +2i \partial_x P_+ \Bigl[W_1
  \partial_x P_-\Bigl(e^{-iF_1/2} \partial_x P_+(e^{iF_1/2}
  - e^{iF_2/2}) \Bigr) \Bigr] \nonumber \\
  & & +
  2i \partial_x P_+ \Bigl[W_1 \partial_x P_-\Bigl((e^{iF_1/2}
  - e^{iF_2/2}) \partial_x P_+(e^{iF_2/2}
   \Bigr)\Bigr] +2i\partial_x P_+ \Bigl(Z \partial_x P_-(e^{-iF_2/2} \overline{w_2} )  \Bigr) \nonumber \\
& & \hspace*{-20mm} +
2i\partial_x P_+ \Bigl[Z \partial_x P_-\Bigl( e^{-iF_2/2} \partial_x P_+
(e^{iF_2/2})\Bigr)\Bigr]
    +\frac{i}{4} \Bigl( P_0(u_1^2) z+ P_0(u_1^2-u_2^2) w_2\Bigr)\; . \label{o14}
 \end{eqnarray*}
 \arraycolsep5pt
 This expression seems somewhat complicated but actually each
  term can be treated  as in Section \ref{Sec4}. We extend the functions
   $ w_i $ and $ F_i $ in the same way as in Section \ref{Sec42}.
   To deal with the difference $ e^{i{\tilde F_1}/2}
  - e^{i{\tilde F_2}/2} $ we use that formally
  $$
  e^{i{\tilde F_1}/2}- e^{i{\tilde F_2}/2} = \sum_{k\in \N} \frac{(i/2)^k}{k!}
   ({\tilde F_1}^k-{\tilde F_2}^k)
   = \sum_{k\ge 1} \frac{(i/2)^k}{k!} ({\tilde F_2}-{\tilde F_2}) \Bigl(
  \sum_{j=0}^{k-1} {\tilde F_1}^{j} {\tilde F_2}^{k-1-j}  \Bigr)
  $$
Moreover, as in \re{lala1} we have
$$
\|P_3({\tilde F_1}-{\tilde F_2})\|_{A_\lambda} \lesssim
 \Bigl\|{\cal F}^{-1}_x \Bigl(\partial_x^{-1}(u_1(0)-u_2(0))\Bigr)\Bigr\|_{L^1_\xi}
  + \|{\tilde {\tilde u}_1}-{\tilde {\tilde u}_2} \|_{N_{\lambda}}( \|{\tilde {\tilde u}_1} \|_{N_{\lambda}}+ \|{\tilde {\tilde u}_2} \|_{N_{\lambda}})
 $$
 and thus
$$
\|{\tilde F_1}-{\tilde F_2}\|_{A_\lambda} \lesssim
 \Bigl\|{\cal F}^{-1}_x \Bigl(\partial_x^{-1}(u_1(0)-u_2(0))\Bigr)\Bigr\|_{L^1_\xi} + \|{\tilde {\tilde u}_1}-{\tilde {\tilde u}_2} \|_{N_\lambda}(1+ \|{\tilde {\tilde u}_1} \|_{N_\lambda}+ \|{\tilde {\tilde u}_2} \|_{N_\lambda}) \; .
$$
Therefore, on account of  Lemmas \ref{line2}-\ref{line3}, \ref{nonlin1}-\ref{nonlin2} and \re{termsup}, we  infer that, for
 $ 0\le s\le 1/2 $,
 \begin{eqnarray*}
\|z\|_{M^{s}_{1,\lambda}} & \lesssim  &
\|z(0)\|_{H^{s}_\lambda}+  e^{{\tilde K}_1+{\tilde K}_2} \Bigl[ \|w_1
\|_{X^{1/2,s}_{1,\lambda}}
 \Bigl( \|z\|_{X^{1/2,0}_{1,\lambda}}
\|u_1\|_{N_{1,\lambda}} \\
& & +(\|\widehat{\partial_x^{-1} v_0}\|_{L^1_\xi}
+ \|v\|_{N_{1,\lambda}}+ \|v\|_{N_1,\lambda}^2)
  (\|w_2\|_{X^{1/2,0}_{1,\lambda}}+
\|u_1\|_{N_1,\lambda}+\|u_2\|_{N_1,\lambda})\Bigr)  \\
& & + \|z\|_{X^{1/2,s}_{1,\lambda}}(\|w_2\|_{X^{1/2,0}_{1,\lambda}}+ \|u_2\|_{N_{1,\lambda}}
+\|u_2\|_{N_{1,\lambda}}^2) +\|z\|_{X^{1/2,0}_{1,\lambda}}\|u_1\|_{N_{1,\lambda}}^2
\\
 & &
+\|v\|_{N_{1,\lambda}}\|u_2\|_{N_{1,\lambda}}\|w_2\|_{X^{1/2,0}_{1,\lambda}}
\Bigr]
\quad .
 \end{eqnarray*}
 where
 $$
 {\tilde K}_1+{\tilde K}_2= C\Bigl(\|\widehat{\partial_x^{-1}u_1(0)}\|_{L^1_\xi}+
\|\widehat{\partial_x^{-1}u_2(0)}\|_{L^1_\xi}+
\|u_1\|^2_{N_{1,\lambda}}+ \|u_2\|^2_{N_{1,\lambda}}\Bigr) \quad .
$$
 Thanks to \re{oo2}  we thus obtain that
 \begin{eqnarray}
 \|z\|_{M^{s}_{1,\lambda}} & \lesssim  & \Bigl(1+\|\varphi_2\|_{L^2_\lambda}+\|\varphi_1\|_{L^2_\lambda}(1+
\lambda^{1/2}) \Bigr)\|\varphi_1-\varphi_2\|_{H^{s}_\lambda}
\nonumber \\
& &  + \varepsilon^2 \, e^{2 C K} \,\Bigl[
 \|w_1\|_{X^{1/2,s}_{1,\lambda}} (\|z\|_{X^{1/2,0}_{1,\lambda}}
 + \|\widehat{\partial_x^{-1} v(0)}\|_{L^1_\xi}
+   \|v\|_{N_{1,\lambda}}) \nonumber \\
 & & + \|z\|_{X^{1/2,s}_{1,\lambda}} + \|v\|_{N_{1,\lambda}}  \Bigr]\quad
 ,\label{Lip1}
 \end{eqnarray}
 since, by Lemma \ref{non1}, it can be easily seen that
\begin{eqnarray*}
\|z(0)\|_{H^{s}_\lambda} & \lesssim &
\|\varphi_1-\varphi_2\|_{H^{s}_\lambda}\Bigl(1+\|\varphi_1\|_{L^2_\lambda}
 +\|\varphi_2\|_{L^2_\lambda}\Bigr)\nonumber\\
 & & +\|e^{-iF_1(0)}-e^{-iF_2(0)}\|_{L^\infty_\lambda} \|\varphi_1\|_{H^{s}_\lambda}
 (1+\|\varphi_1\|_{L^2_\lambda})
\label{o5}
\end{eqnarray*}
with
\begin{eqnarray*}
\|e^{-iF_1(0)}-e^{-iF_2(0)}\|_{L^\infty_{\lambda}}
   \lesssim  \|\partial_x^{-1} (\varphi_1-\varphi_2)\|_{L^\infty_{\lambda}}
   \lesssim
  \lambda^{1/2} \|\varphi_1-\varphi_2 \|_{L^2_\lambda}\quad . \label{o6}
\end{eqnarray*}
On the other hand, proceeding as in Section \ref{Sec5} and using \re{oo2}, one can    check that
\begin{equation}
\|v\|_{N_{1,\lambda}}  \lesssim   \|v(0)\|_{L^2_\lambda}+ \Bigl[
\|z\|_{M^{0}_{1,\lambda}} +\varepsilon^2 \Bigl(
\|\widehat{\partial_x^{-1} v(0)}\|_{L^1_\xi}+
 \|v\|_{N_{1,\lambda}}\Bigr)\Bigr] e^{2CK}\; . \label{Lip2}
\end{equation}
Noticing that by Cauchy-Schwarz in $ \xi $,
 $$
 \|\widehat{\partial_x^{-1} v(0)}\|_{L^1_\xi}\lesssim \lambda^{1/2}
 \|\widehat{ v(0)}\|_{L^2_\xi}\sim \lambda^{1/2}\|v(0)\|_{L^2_\lambda} \
  $$
 and gathering \re{Lip1} and \re{Lip2} we obtain
\begin{equation}
\|v\|_{N_{1,\lambda}}  +\|z\|_{M^0_{1,\lambda}} \lesssim
e^{2CK}(1+\varepsilon^2 \lambda^{1/2})
\|\varphi_1-\varphi_2\|_{L^2_\lambda} \quad . \label{Lip3}
\end{equation}
Coming back to \re{Lip1} this leads to
\begin{equation}
\|z\|_{M^s_{1,\lambda}} \lesssim  e^{2CK}\Bigl(1+\varepsilon^2
\lambda^{1/2} \Bigr)\|\varphi_1-\varphi_2\|_{H^{s}_\lambda} \quad
. \label{Lip4}
\end{equation}
Now, proceeding as in \re{A3}, we infer that
\arraycolsep2pt
\begin{eqnarray*}
v & =& \partial_x F_1 -\partial_x F_2 \nonumber \\
 & = & 2 i e^{iF_1/2} \Bigl[ z +\partial_x P_- \Bigl( e^{-iF_1/2}  -e^{-iF_2/2}
 \Bigr)\Bigr] +2i ( e^{iF_1/2} -e^{iF_2/2}) \Bigl(w_2
  +  \partial_x P_- (e^{-iF_2/2} \Bigr) \label{o71}
\end{eqnarray*}
and thus
 \begin{eqnarray}
P_{>1} v
 & = & 2i P_{>1}(e^{iF_1/2}z)
 +2i P_{>1}\Bigl[  P_{>1}(e^{iF_1/2})\partial_x P_-\Bigl( e^{-iF_1/2} -e^{-iF_2/2} \Bigr)\Bigr] \nonumber \\
& &  \hspace*{-14mm} +2i  P_{>1}\Bigl[( e^{iF_1/2} -e^{iF_2/2})   w_2\Bigr]
  + 2i P_{>1}\Bigl[ P_{>1}( e^{iF_1/2} -e^{iF_2/2}) \partial_x P_- (e^{-iF_2/2} ) \Bigr] \label{o8} \quad .
\end{eqnarray}
\arraycolsep5pt
Therefore, by Lemmas \ref{non1}-\ref{non2}, \re{oo2} and  \re{nonlinear3}
\begin{eqnarray*}
\|J^s_x Q_1 v\|_{L^\infty_1 L^2_\lambda} & \lesssim &
 \Bigl(  \|z\|_{Y^s_{1,\lambda}}+
 \varepsilon^2 ( \|v\|_{L^\infty_1 L^2_\lambda}+\varepsilon^2
 \|\partial_x^{-1}v\|_{L^\infty_1 L^\infty_\lambda})  \Bigr) e^{\tilde K} \\
 & \lesssim &\Bigl(  \|z\|_{Y^s_{1,\lambda}}+
 \varepsilon^2 (1+\varepsilon^2 \lambda^{1/2}) \|v\|_{L^\infty_1 L^2_\lambda}  \Bigr) e^{\tilde K}
\end{eqnarray*}
Since on the other hand (see \re{estP1u}),
\begin{equation}
\|P_1\, v\|_{L^\infty_1 L^2_\lambda} \lesssim
\|\varphi_1-\varphi_2\|_{L^2_\lambda} + \|v\|_{L^4_{1,\lambda}}
 \Bigl( \|u_1\|_{L^4_{1,\lambda}}
+\|u_2\|_{L^4_{1,\lambda}}\Bigr),
\end{equation}
we finally deduce from \re{Lip2}-\re{Lip3} that
\begin{equation}
\|J^s_x v\|_{L^\infty_1 L^2_\lambda}   \lesssim
 e^{4CK} (1+\varepsilon^2 \lambda^{1/2})^2 \|\varphi_1-\varphi_2\|_{H^s_\lambda} \quad .
\label{Lip5}
\end{equation}
\section{Proof of Theorem \ref{main}}
We will first prove the local well-posedness result for small data,
the result for arbitrary large data  will
  then follow from scaling arguments.
 \subsection{Well-posedness for  small initial data}
 For any  $K\ge 1 $ and $\lambda\ge 1 $ given,
  let $ u_0\in  B_{K,\lambda}\cap H^s_\lambda $
  with $ 0\le s \le 1/2  $ and
 let  $\{ u_{0}^n \} \subset {H}_0^{\infty}(\T)\cap B_{K,\lambda} $  converging to $ u_0 $ in
 $ H^{s}(\T) $. We denote by $ u_n $ the solution of (BO) emanating from $ u_{0}^n $.
 From standard  existence theorems (see for instance \cite{ABFS}, \cite{Io}),  $ u_n\in C(\R;{H}^\infty_{0,\lambda}) $.
  According to
 \re{o2} and \re{ets}, for all $n\in \N^* $,
 \begin{equation}
   \|u_n \|_{N_{1,\lambda}}
    +  \|w_n \|_{M^{0}_{1,\lambda}}\lesssim
  e^{2CK} \,  \varepsilon(K)^{2} \label{bound1}
   \end{equation}
   and
   \begin{equation}
  \|u_n\|_{L^\infty_1 H^s_\lambda}
    +  \|w_n \|_{M^{s}_{1,\lambda}}
    \lesssim  e^{2C\|\widehat{\partial_x^{-1} u_0}\|_{L^1_\xi}}
    \, \|u_0\|_{H^s_\lambda}
    \; , \label{bound2}
    \end{equation}
    where $ w_n =\partial_x P_+(e^{-iF_n/2}) $
     is the gauge transform of $w_n $.
   Note that  this uniform bound would enable
      to prove the local existence
  for $ s>0 $ by using weak convergences. On the other hand,  for $s=0 $,
   weak convergences  would   not be sufficient to pass to the limit
    on the nonlinear term $ u^2 $. Actually,
 with \re{Lip3} and \re{Lip5} in hand, we observe that the approximative sequence
$ u^n $ constructed for the local existence result is a Cauchy sequence
 in  $ C([0,1]; {H}^{s}_{0,\lambda}) \cap N_{1,\lambda} $
 since the $u_n $ satisty \re{o2}-\re{ets} and
 $ u_{0,n} $ converges to $ u_0 $ in $ {H}^{s}_{0,\lambda}$.
  Hence, $ u_n $ converges strongly to some $ u $ in
  $ C([0,1]; {H}^{s}_{0,\lambda})\cap N_{1,\lambda} $. This strong
  convergence permits to pass easily to the limit on the nonlinear term and thus $ u $ is  a solution of (BO). Moreover, from \re{ets} and \re{Lip4} it follows that the sequence of gauge transforms $ w_n $ of $ u_n $ is a Cauchy sequence in $ M^s_{1,\lambda} $. Hence
   $ w_n= \partial_x P_+(e^{-iF_n/2}) $ converges toward some function $ w $ in $ Y^s_{1,\lambda} $ and from the strong convergence of $ u $ it is easy to check that $ w= P_+(e^{-iF/2}) $
    with $ F=\partial_x^{-1} u $.
   \\
  Now let $ u^1$ and $ u^2 $ be two solutions emanating from $ u_0 $ belonging
   to $ N_{1,\lambda}$ such that their associated  gauge functions
    belong to $ X^{1/2,0}_\lambda $. According to \re{nonlinear1}, the gauge functions
     belong in fact to $ M^0_{1,\lambda} $ and using the same dilation argument we used to prove
    the uniform boundness of the solution, we can show that for $\beta $ large enough
     and $i=1,2 $,
     $$
   \| u_\beta^i \|_{N^{0}_{1,\lambda\beta}}
     + \|  w_\beta^i \|_{M^{0}_{1,\lambda\beta }} \lesssim e^{2CK} \,
    \|u_{0,\beta}\|_{L^2_{\lambda\beta} }\lesssim e^{2CK} \beta^{-1/2} \varepsilon(K)^2 \quad
  $$
  with $K =\|\widehat{\partial_x^{-1} u_{0,\beta}}\|_{L^1_\xi}+1=\|\widehat{\partial_x^{-1}
   u_0}\|_{L^1_\xi}+1$.
 Therefore, for $ \beta $ large enough, $( u^i_{\beta}, w^i_{\beta}) $ satisfies the smallness condition \re{oo2} with
  $\varepsilon=\varepsilon(K)$ and  $ u_{0\beta}\in B_{K,\lambda\beta} $.
 It then follows from \re{Lip3} that  $ u_\beta^1\equiv u_\beta^2 $ on $ [0,1] $
  and thus $ u^1\equiv u^2 $ on $ [0,1/\beta^2] $. This proves the uniqueness result for  initial data belonging to
   $ B_{K,\lambda}$.
Moreover, \re{Lip5} clearly ensures
 that  the  flow-map is
Lipschitz from  $ B_{K,\lambda}\cap H^s_\lambda $ into
  $ C([0,1];{H}^{s}_{0,\lambda}) $.
  \subsection{The case of arbitrary large initial data}
We use again  the dilation invariance of (BO) to extend the result for arbitrary large data.
Recall that  that if $ u(t,x) $ is a $2\pi $-periodic solution of (BO) on $[0,T] $ with initial data
 $ u_0 $ then $ u_\lambda(t,x)=\lambda^{-1} u(\lambda^{-2}t,\lambda^{-1} x)$ is a $ (2\pi \lambda) $-periodic solution of (BO) on
 $ [0,\lambda^2 T] $
 emanating from $u_{0,\lambda}=\lambda^{-1} u_0(\lambda^{-1} x) $. Recall also
 that the associated
  gauge functions satisfy $ w_\lambda(t,x)=\lambda^{-1} w(\lambda^{-2}t,\lambda^{-1} x) $.\\
Let $ u_0\in {H}_0^{s}(\T) $ with  $ 0\le s\le 1/2 $. Note that
$$
 \|\widehat{\partial_x^{-1} u_0} \|_{L^1_\xi} \lesssim \|u_0\|_{L^2_1} \quad .
 $$
 We thus set $ K=\|u_0\|_{L^2_1} +1$ and  take
$$
\lambda = \max\Bigl(1,\varepsilon(K)^{-4} \| u_0\|_{L^2_1}^2
\Bigr) \ge 1
$$
 so that
 $$
 \|u_{0,\lambda}\|_{L^2_\lambda} \le \lambda^{-1/2} \|u_0\|_{L^2_1} \le \varepsilon(K)^2 \; .
 $$
Recalling that $ \|\widehat{\partial_x^{-1} u_{0,\lambda}}
\|_{L^1_\xi}= \|\widehat{\partial_x^{-1} u_0} \|_{L^1_\xi}$, it
follows that $ u_{0,\lambda} $ belongs to $ B_{K,\lambda} $ and
so we  are reduced to the case of small initial data. Therefore,
there exists a unique
 solution $ u_\lambda\in C([0,1];{H}^{s}_{0,\lambda})\cap N_{1,\lambda} $ of (BO)
  with $ w_\lambda\in M^s_{1,\lambda} $. This proves the
existence and uniqueness of
  the solution $ u $ of (BO) in the class
 $$
 u\in  C([0,T];{H}_0^{s}(\T))\cap N_{T,1}
  ,  \quad w \in M^{s}_{T,1}
  $$
  emanating from $ u_0 $ where $ T=T(\|u_0\|_{L^2})$ and $ \alpha\mapsto T(\alpha)  $ is a non increasing
   function on $\R_+^* $.
   The fact that the flow-map is Lipschitz on every bounded set of $ {H}_0^{s}(\T) $ follows as well since $\lambda $  only depends on $ \|u_0\|_{L^2}$.\vspace*{3mm}\\

 Note that the change of unknown \re{chgtvar} preserves the continuity of the solution and
 the continuity of the flow-map in $ H^{s}(\T)$. Moreover, the Lipschitz property
 (on bounded sets) of the flow-map is also preserved on the hyperplans of $ H^{s}(\T) $
  of functions
   with fixed  mean value. Finally, the global well-posedness result follows directly by
   combining the conservation of the $ L^2 $-norm  and   the local
    well-posedness result.

 \section{Proof of Theorem \ref{illposedness}}
 \subsection{Analycity of the flow-map}
Let us   prove the analyticity of the solution-map $ \Psi \, :\, u_0 \mapsto u $
 from $ { H}_0^s(\T) $ to $ C([0,1];H^s(\T)) $ at the origin.
  Note that the other points of $ { H}_0^s(\T) $
   could be handle
 in the same way. Also we restrict ourself to the case $ 0\le s\le 1/2 $ but the case $ s\ge 1/2 $ can
  be treated in a similar way (in fact easier) by using the results of \cite{M1}. \\
  The analycity of the flow-map will be a direct consequence
   of the three following ingredients :\\
  $\bullet $ The Lipschitz property of $\Psi $ proven in Section \ref{section6}.\\
  $\bullet$  The fact that it appears
  only polynomial or analytic functions of $ u $ in the equations we
   deal with. \\
  $\bullet $ We have an absolute convergence, in the norms we are interested in, of the serie
   obtained by replacing the analytic functions of $ u $ by their associated entire series.
   \vspace*{2mm} \\
  So, let $\varphi\in { H}_0^s(\T) $ with $ \| \varphi\|_{H^s_1}=1 $ and let $ \varepsilon>0 $ be a small
   real number to be fixed later. Taking $ u_0=\varepsilon \varphi $ we know from \re{Lip3}, \re{Lip4} and \re{Lip5} that,
    for $\varepsilon $ small enough, there exists $ c_1>0 $ such
    that
    the corresponding solution $ u $ and its gauge transform $ w $ verify
    \begin{equation}
    \|u \|_{N_{1,1}}+\|u\|_{L^\infty_1 {H}^s_1}+ \|w\|_{M^s_{1,1}}
    \le  c_1 \, \varepsilon  \quad ,\label{huhu1}
    \end{equation}
  Now let  $ C>0 $ be  a universal constant we take  very large
     (We can take for example
      $ C >0 $ to be  the exponential of the sum of all
     the implicit constants interfering in our estimates in Sections \ref{Sec4}-\ref{Sec5}). According to \re{bobo} and \re{otot}, we  get
     $$
\|P_{3} u-\varepsilon V(t) P_{3} \varphi\|_{N_{1,1}} \le
  C (c_1 \, \varepsilon)^2 \quad.
     $$
    On the other hand, since $\partial_x^{-1}\varphi $ belongs to $ H^{s+1}$ which is an algebra, it holds
     in $ H^{s+1}_1 $
    $$
    W(0)= P_+(e^{-i\varepsilon\partial_x^{-1} \varphi/2}) =1-\frac{i}{2}\varepsilon
    P_+(\partial_x^{-1} \varphi) + \sum_{k\ge 2}
     (\frac{-i\varepsilon}{2})^{k} \frac{1}{k!} P_+\Bigl((\partial_x^{-1} \varphi)^k \Bigr) \quad .
      $$
 and thus
 $$
 w(0)= -\frac{i}{2}\varepsilon P_+ (\varphi) +\Lambda_\varepsilon \mbox{ with }
  \|\Lambda_\varepsilon\|_{H^s_1} \le  4 \varepsilon \quad .
 $$
 Consequently,
 $$
 V(t) w(0) = -\frac{i}{2}\varepsilon V(t) P_+ \varphi +V(t) \Lambda_\varepsilon \mbox{ with }
  \|V(t)\Lambda_\varepsilon\|_{M^s_{1,1}}
  \le    C(4  \varepsilon)^2\quad .
  $$
  Now according to \re{C1}, \re{mi2}, \re{mi11} and \re{termsup},
   we infer that
  $
  \| w-V(t) w(0) \|_{M^s_{1,1}} \le\, C    (c_1
  \varepsilon)^2
  $
  and thus
 \begin{equation}
  \| w+\frac{i}{2}\varepsilon V(t) P_+ \varphi\|_{M^s_{1,1}} \le 2 C
    (c_1 \varepsilon)^2 \quad . \label{huhu2}
 \end{equation}
 It then follows from \re{A4}-\re{mi61}, \re{defF}, \re{co3} and
  \re{Xnorm} that
 $$
 P_{>3}(u)=2i P_{>3}w +\tilde{\Lambda}_\varepsilon=\varepsilon V(t)
 P_{>3} (\varphi)
  + \tilde{\tilde \Lambda}_\varepsilon
  $$
  for some function $ \tilde{\tilde \Lambda}_\varepsilon $ satisfying
   $ \|\tilde{\tilde \Lambda}_\varepsilon\|_{N_{1,1}}
   +\|\tilde{\tilde \Lambda}_\varepsilon\|_{L^\infty_1 H^s}
   \le 3C(c_1  \varepsilon)^2 $.\\
   We thus finally get,
   \begin{equation}
  \|u-\varepsilon V(t) \varphi\|_{N_{1,1}}+
    \|u-\varepsilon V(t) \varphi\|_{L^\infty_1 H^s_1}+
    \| w+\frac{i}{2}\varepsilon V(t) P_{>3} \varphi \|_{M^s_{1,1}} \le
    6C(c_1 \varepsilon)^2 \; . \label{huhu3}
 \end{equation}
 In the same way, according to \re{eq2w}, expanding $ e^{-iF/2} $ and $ e^{iF/2} $ as in
 Section \ref{Sec4}, with
 \re{huhu1}-\re{huhu3}  in hand, we get
 \begin{eqnarray*}
 w & = & -\frac{i}{2}\varepsilon V(t) P_+(\varphi) -\varepsilon^2\Bigl[
  \frac{1}{4} V(t) P_+(\varphi \partial_x^{-1} \varphi) + 2i \int_0^t V(t-t') \partial_x P_+
  \Bigl(W_1 \partial_x P_-(\overline{w_1})\Bigr) \\
  & & +\Lambda_\varepsilon \; .
 \end{eqnarray*}
 where
 $$
 u_1=V(t) \varphi, \; W_1=  -\frac{i}{2} V(t) P_+ (\partial_x^{-1} \varphi), \;
 w_1= \partial_x W_1
 $$
 and  $ \|\Lambda_\varepsilon\|_{M^s_1} \lesssim 6C^2 \,
 (c_1\varepsilon)^3 $
   and so on ...\\
 Iterating this process we obtain that there exists $\varepsilon_0>0 $ such that the following
   asymptotic expansion of $ u$ in term of
 $\varphi $ holds  absolutely in $ C([0,1];H^s(\T)) $
  for $ 0<\varepsilon\le \varepsilon_0 $,
 \begin{equation}
 u= \sum_{k\ge 1} \varepsilon^k A_k(\varphi) \label{deve} \quad .
 \end{equation}
 Here, $ A_1(\varphi)= t\mapsto V(t) \varphi $ and more generaly $ A_k $ is a
 continuous k-linear operator from $ { H}_0^s(\T) $ to $ C([0,1];H_0^s(\T)) $.
 Therefore $ u $ is real-analytic and in particular $ C^\infty $
  at the origin of $ { H}_0^s(\T) $.
  Moreover, since
  $$
  u(t,\cdot)= \varepsilon U(t) \varphi+\frac{1}{2} \int_0^t V(t-t') \partial_x u^2(t')\, dt'\; ,
  $$
  by identification we infer that
  \begin{equation}
  A_k(\varphi)=t\mapsto\frac{1}{2} \sum_{k_1,k_2\ge 1 \atop k_1+k_2=k}
  \int_0^t V(t-t') \partial_x \Bigl(A_{k_1}(\varphi) A_{k_2} (\varphi)\Bigr)(t') \, dt'
  \label{kitere}
  \end{equation}
  \subsection{Non smoothness of the flow-map in $ H^s(\T) $, $ s<0 $.}
  Let us start by computing  $ A_k(t,\lambda \cos(Nx)) $ for  k=1,2,3. Of course,
$$
A_1(t,cos(Nx))=\cos(Nx-N^2t) \quad .
$$
Since $ \partial_x \Bigl(A_1(t,cos(Nx))\Bigr)^2= -N \sin(2Nx-2N^2t) $ we infer
that
\begin{eqnarray*}
A_2(t,\cos(Nx)) & =& \frac{1}{2}\int_0^t V(t-t')
\partial_x \Bigl(A_1(t,cos(Nx))\Bigr)^2 (t') \, dt \\
 & = & -\frac{N}{2} \int_0^t \sin\Bigl(2Nx-2N^2t'-4N^2(t-t')\Bigr)\, dt' \\
  &=& \frac{1}{4N} \Bigl[ \cos(2Nx-2N^2t)-\cos(2Nx-4N^2t) \Bigr]
\end{eqnarray*}
In the same way,
\begin{eqnarray*}
\partial_x\Bigl( A_1(1,\cos(Nx)) A_2(t,\cos(Nx)) \Bigr) & =&
   -\frac{1}{8} \Bigl[ \sin(Nx-N^2t)-\sin(Nx-3N^2t) \Bigr] \\
   & & -\frac{3}{8} \Bigl[ \sin(3Nx-3N^2t)-\sin(3Nx-5N^2t) \Bigr]
\end{eqnarray*}
and thus
\begin{eqnarray*}
A_3(t,\cos(Nx)) & =& \int_0^t V(t-t')
\partial_x \Bigl(A_1(t,cos(Nx))A_2(t,\cos(Nx))\Bigr) (t') \, dt \\
 & = & -\frac{1}{8} \int_0^t
 \Bigl[ \sin(Nx-N^2t)-\sin(Nx-3N^2t'-N^2(t-t')) \Bigr]\, dt'\\
&  & -\frac{3}{8} \int_0^t
 \Bigl[ \sin(3Nx-3N^2t'-9N^2(t-t'))-\sin(Nx-5N^2t'-N^2(t-t')) \Bigr]\, dt'\\
& =& -\frac{t}{8}  \sin (Nx-N^2t) \\
&  & +\frac{1}{16N^2} \Bigl[ \cos(Nx-3N^2t)-\cos(Nx-N^2t) \Bigr]\\
 & & +\frac{1}{16N^2} \Bigl[ \cos(3Nx-3N^2t)-\cos(3Nx-9N^2t) \Bigr] \\
 & & -\frac{3}{32 N^2} \Bigl[ \cos(3Nx-5N^2t)-\cos(3Nx-9N^2t) \Bigr]
\end{eqnarray*}
Therefore, setting $ \Psi_N=N^{-s} \cos(Nx)  $ it follows that
$$
\|A_3(t,\Psi_N)\|_{H^s}  \gtrsim t\, N^{-2s} \|\Psi_N\|_{H^s}^3
$$
 and from
 standard considerations (cf. \cite{Bo2})   the flow-map cannot be of class $C^3 $
at the origin from $ { H}_0^s(\T) $ into $ { H}_0^s(T) $ as
soon as $ s< 0$. Moreover, by a direct induction argument it is not too hard to check that for any $ k\ge 4 $,
$$
\|A_k(t,\cos(Nx))\|_{H^s} \le  {\tilde C}_k \, N^{s} \quad .
$$
Therefore, for any fixed integer $ K\ge 4 $,
$$
\Bigl\|\sum_{k=4}^{K+2} A_k(t,\varepsilon\cos(Nx)) \Bigr\|_{H^s} \le C_K \varepsilon^4 N^ {s} \quad .
$$
 Now, taking as initial data $ \varphi_N=\varepsilon_N \cos(Nx) $ with   $ 0<\varepsilon_N\le \varepsilon_0/2 $,
   we know from \re{deve} that the associated solution $ u_N $ can be written in $ L^2(\T) $  as
   $$
   u_N(t,\cdot)=\sum_{k\ge 1} \varepsilon_N^k A_k(t,\cos(Nx)) \quad .
   $$
  For $ N $ large enough
and $ s\le 0 $, we thus  deduce from the computation of $ A_2(t,\cos(Nx)) $
 and  $ A_3(t,\cos(Nx)) $ above that
\begin{eqnarray*}
\|u_N(t,\cdot) -V(t)\varphi_N \|_{H^s} & \gtrsim & t \varepsilon_N^3
\|
\sin(N x-N^2 t) \|_{H^s} -2 \varepsilon_N^2 N^{s-1}-C_K \varepsilon_N^4 N^ {s} \\
& & -{\tilde C}\sum_{k=K+3}^\infty (\frac{\varepsilon_N}{\varepsilon_0})^{k} \|A_k(t,\varepsilon_0\cos(Nx))\|_{L^2} \\
& \gtrsim & t \varepsilon_N^3
\|\sin(N x-N^2 t) \|_{H^s} -2 \varepsilon_N^2 N^{s-1}-C_K \varepsilon^4_N N^ {s} -C \varepsilon_N^{K+3}\\
& \gtrsim &  \varepsilon^3_N N^{s} \Bigl( t-\frac{2}{N\varepsilon_N}
-C_K \varepsilon_N -C \varepsilon_N^{K}N^{-s} \Bigr) \quad .
\end{eqnarray*}
For any  $0<\alpha <1 $ and $ s<0 $ fixed, we take $ K>0 $ such that
$$
\frac{|s|}{K} <1 \quad \mbox{ and } \quad \frac{4}{K} <\alpha
\quad .
$$
Setting
$$
\varepsilon_N= \min\Bigl( \frac{\varepsilon_0}{2}, \frac{t}{4C_K},\,  (\frac{t\, N^s}{4C})^\frac{1}{K} \Bigr)
$$
we infer that for $ N $ large enough, 
\begin{eqnarray*}
\|u_N(t,\cdot) -V(t)\varphi_N \|_{H^s}  & \gtrsim &  t \varepsilon^3_N N^{s} \\
& \gtrsim & t \varepsilon_N^{2-\alpha} N^{-\alpha  s} \|\varphi_N \|_{H^s}^{1+\alpha}\\
& \gtrsim & t N^{-\frac{\alpha s}{2}}  \|\varphi_N \|_{H^s}^{1+\alpha} \quad .
\end{eqnarray*}
It follows that the flow-map (if it coincides with the standard flow-map on $ { H}_0^\infty(\T)$) cannot be of class $
C^{1+\alpha} $ at the origin from $ { H}_0^s(\T) $ into $ {
H}_0^s(\T) $.
\section{Appendix}
\subsection{Proof of Lemma \ref{non1}}
We separate the low and the high modes of $ h $. To treat the high modes part,
 we observe that by  Leibniz
rule for fractional derivatives (cf. \cite{KPV4}) and Sobolev inequality,
\begin{eqnarray*}
\Bigl\|J_x^\alpha \Bigl( Q_1( h) \, g \Bigr) \Bigr\|_{L^q_\lambda} & \lesssim &
 \Bigl\|J_x^\alpha  Q_1 (h ) \Bigr\|_{L^{4/\alpha}_\lambda}
  \|g \|_{L^\frac{4q}{4-\alpha q}_\lambda} + \|h\|_{L^\infty_\lambda} \|J_x^\alpha g \|_{L^q_\lambda} \\
  & \lesssim &
 \Bigl\|J_x^{\alpha+ 1/2-\alpha/4} Q_1 (h ) \Bigr\|_{L^{2}_\lambda}
  \|J_x^{\alpha/4} g \|_{L^q_\lambda} + \|h\|_{L^\infty_\lambda} \|J_x^\alpha g \|_{L^q_\lambda} \\
  & \lesssim & ( \|\partial_x h \|_{L^2_\lambda} +\|h\|_{L^\infty_\lambda} )
  \|J_x^\alpha g \|_{L^q_\lambda}\; .
\end{eqnarray*}
On the other hand, one can easily check that
$$
\Bigl\|J_x^1 \Bigl( P_1 (h) \, g \Bigr) \Bigr\|_{L^q_\lambda} \lesssim
 (\|h\|_{L^\infty_\lambda} +\|\partial_x h \|_{L^2_\lambda})
  \|J_x^1 g \|_{L^q_\lambda} \quad \mbox{ and } \quad
  \| P_1(h) g \|_{L^q_\lambda} \lesssim \|h\|_{L^\infty_\lambda}\| g \|_{L^q_\lambda}\; .
$$
Interpolating between this two estimates we obtain the desired
 estimate on the low modes part.
 \subsection{Proof of Lemma \ref{non3}} \label{proofnon2}
Clearly the low modes part  of $ z v $ can be estimated directly by an Holder inequality.
  Now, using the nonhomogeneous Littlewood-Paley decomposition, we get for $ q\ge 8 $,
\begin{eqnarray*}
\Delta_q(zv) &  = &  \sum_{|i|\le 2}   \Delta_q \Bigl(
\Delta_{q-i}(v) \sum_{j=0}^{q-i-2} \Delta_j(z) \Bigr) \\
& & \hspace*{-18mm}+
 \sum_{|i|\le 2}   \Delta_q \Bigl(
\Delta_{q-i}(z) \sum_{j=0}^{q-i-2} \Delta_j(v) \Bigr)
  +  \Delta_q  \Bigl( \sum_{i\ge q-2} \sum_{|j|\le 1} \Delta_{i-j}(v) \Delta_i(z)\Bigr)  \quad .
 \end{eqnarray*}
 Therefore,
 \begin{eqnarray}
\sum_{q\ge 8} \|\Delta_q(zv)\|_{L^4_{t,\lambda}}^2  & \lesssim &
 \|z\|_{L^\infty_{t,\lambda}}^2 \sum_{q\ge 4} \|\Delta_q(v)\|_{L^4_{t,\lambda}}^2
 \nonumber \\
 & &\hspace*{-18mm} + \|v\|_{L^4_{t,\lambda}} \Bigl( \sum_{q\ge 4}
  \|\Delta_q(z)\|_{L^\infty_{t,\lambda}}^2
  + \sum_{q\ge 4} \sum_{k\ge q-2}
   \|\Delta_k(z)\|_{L^\infty_{t,\lambda}}^2 \Bigr) \; .\label{yoyo}
 \end{eqnarray}
 The  desired result follows since for $ k \ge 2 $,
 $$
\|\Delta_k(z)\|_{L^\infty_{t,\lambda}} \lesssim
 2^{-k/4} \| z_x\|_{L^\infty_t  L^2_\lambda} \;.
$$

\vskip0.3cm \noindent{\bf Acknowledgments .} The author is very grateful to Nikolay Tzvetkov for fruitful discussions on the problem. He is also grateful to the Referees for several valuable suggestions.

\end{document}